\documentclass[preprint,11pt]{imsart}
\RequirePackage{amsthm,amsmath,
	amssymb,bbm}
\RequirePackage[numbers]{natbib}%numbers
\RequirePackage[colorlinks,citecolor=blue,urlcolor=blue]{hyperref}
\usepackage{geometry}
\startlocaldefs
\numberwithin{equation}{section}
\theoremstyle{plain}
\newtheorem{theorem}{Theorem}[section]
\newtheorem{assume}{Assumption}
\newtheorem{lemma}{Lemma}
\newtheorem{cor}{Corollary}[section]

\newtheorem{dfn}{Definition}
\newtheorem{rmk}{Remark}

\endlocaldefs
\allowdisplaybreaks
\begin{document}
	\begin{frontmatter}
		\title{Adaptive posterior concentration rates for sparse high-dimensional linear regression with random design and unknown error variance}
		%\title{A sample article title with some additional note\thanksref{t1}}
		%\runtitle{sparse high dimensional linear regression with random design and unknown variance error}
		%\thankstext{T1}{A sample additional note to the title.}
		
\begin{aug}
	%%%%%%%%%%%%%%%%%%%%%%%%%%%%%%%%%%%%%%%%%%%%%%%
	%% Only one address is permitted per author. %%
	%% Only division, organization and e-mail is %%
	%% included in the address.                  %%
	%% Additional information can be included in %%
	\author{\fnms{The Tien}~\snm{Mai}\ead[label=e1]{the.t.mai@ntnu.no}\orcid{0000-0002-3514-9636}}	%%%%%%%%%%%%%%%%%%%%%%%%%%%%%%%%%%%%%%%%%%%%%%
	%% Addresses                                %%
	%%%%%%%%%%%%%%%%%%%%%%%%%%%%%%%%%%%%%%%%%%%%%%
\address{
	Department of Mathematical Sciences, Norwegian University of Science and Technology,
	\\
	Trondheim 7034, Norway.
	\\
	\printead[presep={\ }]{e1}
}
\runauthor{Mai T. T.}
\end{aug}
\begin{abstract}	
This paper investigates sparse high-dimensional linear regression, particularly examining the properties of the posterior under conditions of random design and unknown error variance. We provide consistency results for the posterior and analyze its concentration rates, demonstrating adaptiveness to the unknown sparsity level of the regression coefficient vector. Furthermore, we extend our investigation to establish concentration outcomes for parameter estimation using specific distance measures. These findings are in line with recent discoveries in frequentist studies. Additionally, by employing techniques to address model misspecification through a fractional posterior, we broaden our analysis through oracle inequalities to encompass the critical aspect of model misspecification for the regular posterior. Our novel findings are demonstrated using two different types of sparsity priors: a shrinkage prior and a spike-and-slab prior.
\end{abstract}
\end{frontmatter}
	%%%%%%%%%%%%%%%%%%%%%%%%%%%%%%%%%%%%%%%%%%%%%%
	%% Please use \tableofcontents for articles %%
	%% with 50 pages and more                   %%
	%%%%%%%%%%%%%%%%%%%%%%%%%%%%%%%%%%%%%%%%%%%%%%
	\tableofcontents
\section{Introduction}
In this paper, we study the following Gaussian linear regression model
\begin{equation*}
y_i = x_i^\top \theta + \sigma \epsilon_i, \quad i=1,\ldots,n, \quad \text{independent}, 
\end{equation*}
where $y_i$ is the response variable, $x_i=(x_{i1},\ldots,x_{ip})^\top \in \mathbb{R}^d $ is a vector of predictor variables,  and $\epsilon_i \sim N(0,1)$ is the random error term. Here, $\theta \in \mathbb{R}^d $ is an unknown vector of regression coefficients and $\sigma > 0$ is an unknown scale parameter. This work focuses on high-dimensional sparse regression, where \( n \ll d \) and only a small subset of the coefficients in the true \( \theta \) are nonzero.

To address challenges in high-dimensional settings, numerous methodologies for parameter estimation have been proposed in the literature, predominantly employing penalization techniques. Among these methodologies, the Lasso stands out as one of the most widely recognized approaches \citep{tibshirani1996regression}. Some exemplary texts on frequentist literature concerning high-dimensional models include: \cite{hastie2009elements,buhlmann2011statistics,giraud2021introduction,wainwright2019high}. More recently, in a notable work by \cite{bellec2018slope}, the authors demonstrated that an adaptive version of Lasso achieves the minimax-optimal rates of convergence of order \( s\log(d/s)/n \) with random design, where \( s \) is the number of nonzero coefficients in the true \( \theta \). Extending this work to the case of unknown error variance \( \sigma \), \cite{derumigny2018improved} showed that an adaptive version of square-root Lasso can attain the same rates.

In high-dimensional linear regression, Bayesian methods have seen a great interest in recent years, for example \citep{jiang2007bayesian,castillo2015bayesian,rovckova2018spike,belitser2020empirical,gao2020general,ning2020bayesian,jeong2021unified,bai2022spike,ray2022variational} among others. These methods offer a flexible framework conducive to incorporating hidden lower-dimensional structures such as sparsity through prior distributions.  More specifically, a notable paper by \cite{castillo2015bayesian} explores the connection between regularization and Bayesian methods, demonstrating that the posterior contraction rates are of order \( s\log(d)/n \). This rate is further improved in a recent work by \cite{gao2020general} to \( s\log(d/s)/n \). For a comprehensive overview of Bayesian methods in this context, readers may refer to the recent review \cite{banerjee2021bayesian}.

However, from a theoretical perspective, these studies predominantly focus on either a fixed random design or known error variance, but not both simultaneously. For instance, \cite{fang2024high} addresses scenarios with unknown error variance but a fixed design, while \cite{jiang2007bayesian} tackles cases involving random design but known error variance. Our objective in this work is to study the behavior of the posterior under both random design and unknown error variance, comparable to the frequentist approach as in \cite{derumigny2018improved}.

In this study, we first focus on the properties of the fractional posterior, also referred to as tempered posteriors, where a fractional power of the likelihood is used, as detailed in \cite{bhattacharya2016bayesian, alquier2020concentration}. It is worth mentioning that using fractional posteriors or replacing the likelihood by a loss function, has received significant attention in recent years, as in \cite{matsubara2022robust, hammer2023approximate, jewson2022general, yonekura2023adaptation, medina2022robustness, mai2017pseudo, grunwald2017inconsistency, bissiri2013general, yang2020alpha, lyddon2019general, syring2019calibrating, Knoblauch, mai2023reduced, hong2020model}. We also note that several works, such as \cite{martin2017empirical, liu2021bayesian, martin2020empirical, fang2024high}, have focused on the fractional posterior in linear models, but only with fixed design.

In addition to our novel focus on random design and unknown error variance, we present novel results concerning the regular posterior in the context of a model misspecification by employing techniques from the fractional posterior. The main results are obtained in the form of oracle inequalities accounting for the model misspecification error. Our novel findings are demonstrated using two different types of sparsity priors: a shrinkage prior and a spike-and-slab prior. We demonstrate that both the regular posterior and the fractional posterior adaptively concentrate at the optimal rate in high-dimensional sparse linear regression without requiring knowledge of the underlying sparsity of the true coefficients. Our results are comparable to those in the frequentist literature \cite{bellec2018slope,derumigny2018improved}, but unlike their work, we do not assume a Restricted Eigenvalue condition.

The rest of the paper is structured as follows. Section \ref{sc_modelmethod} introduces the model, methodology, and priors. Section \ref{sc_concentration} presents consistency results and concentration rates for the fractional posterior. Subsequently, Section \ref{sc_regular_poster_n_horseshoe} provides the consistency results, concentration rates, and properties of the regular posterior in the context of model misspecification. Section \ref{sc_regular_poster_n_horseshoe} also contains results regarding a spike and slab prior. All technical proofs are deferred to Appendix \ref{sc_proofs}. The paper concludes with discussions in Section \ref{sc_discuconcluse}.

\section{Model and method}
\label{sc_modelmethod}
\subsection*{Notations}
For an $\alpha\in(0,1)$, let $P $ and $ R $ be two probability distributions. Let $\mu$ be any measure such that $P\ll \mu$ and $R\ll \mu$.  Defined by
$
D_{\alpha}(P,R)  =
\frac{1}{\alpha-1} \log \int \left(\frac{{\rm d}P}{{\rm d}\mu}\right)^\alpha \left(\frac{{\rm d}R}{{\rm d}\mu}\right)^{1-\alpha} {\rm d}\mu  
$, the $\alpha$-R\'enyi divergence between $P $ and $R $. The Kullback-Leibler divergence is defined by
$
\mathcal{K}(P,R)  = 
\int \log \left(\frac{{\rm d}P}{{\rm d}R} \right){\rm d}P $  if  $ P \ll R
$, and  $
+ \infty $  otherwise. Let $ \|\cdot\|_q $ denote the $ \ell_q $-norm, $ \|\cdot\|_\infty $ denote the max-norm of vectors, and let $ \|\cdot\|_0 $ denote the $ \ell_0 $ (quasi)-norm (the number nonzero entries) of vectors.

\subsection{Model}

In this work,  a linear regression  model is considered. More specially, let $ Z_i:= (Y_i,X_i)\in \mathbb{R} \times \mathbb{R}^d $, $ i=1,\ldots,n $, be a collection of $n$ independent and identically distributed (i.i.d) random samples drawn from a distribution $ P_{\theta,\sigma}, \theta \in \Theta \subset \mathbb{R}^d, \sigma \in  \mathbb{R}^+ $. Under the model $ P_{\theta,\sigma} $, the linear relationship between the response $ Y_i $ and predictor $ X_i $ is such that
\begin{equation}
\label{eq_logistic_model}
Y_i
=
X_i^\top \theta +  \epsilon_i
,
\end{equation}
where $ \epsilon_i \sim N(0, \sigma^2) $ is random noise variable, with unknown $ \sigma >0 $. 

We consider the scenario of random design, assuming that the distribution of \( X_1^n \) is independent of the parameters and that \( X_1 \) is independent of \( \epsilon_1 \).

First, it is assumed that the model is well-specified, meaning there exists a parameter value \( \theta_0 \in \Theta \) and \( \sigma_0 \in \mathbb{R}^+ \) such that \( (Y_i, X_i)_{i=1}^n \) are distributed according to \( P_{\theta_0, \sigma_0}^{\otimes n} \). Specific notations will be provided later if this assumption is relaxed.

Let 
$ L_n(\theta,\sigma ) 
= 
\prod_{i=1}^n p_{\theta,\sigma}(Y_i,X_i) $ be the likelihood and $\pi $ be a prior distribution for $ \theta $ and $ \sigma $. The fractional posterior is given by
\begin{equation}
\label{eq_fractional_posterior}
\Pi_{n,\alpha}(\theta,\sigma )
\propto 
L_n^{\alpha}(\theta,\sigma) \pi(\theta,\sigma)
,
\end{equation}
with $ \alpha \in (0,1) $, as in \cite{bhattacharya2016bayesian,alquier2020concentration}. When $ \alpha =1 $ as in Section \ref{sc_regular_posterior}, we recover the usual posterior distribution. Here we consider that the priors for $ \theta $ and $ \sigma $ are independent and that $ \pi(\theta,\sigma) = \pi_\theta(\theta) \pi_\sigma(\sigma) $.

A high-dimensional sparse scenario is studied in this paper, assuming that $ s^* < n < d $, where $ s^*:= \|\theta_0 \|_0 $.

\subsection{Prior specification}

Choosing an appropriate prior distribution is crucial for attaining a desirable posterior concentration rate in high-dimensional models. In Section \ref{sc_spike_n_slab_prior}, we demonstrate that a spike-and-slab prior yields results comparable to those of the scaled Student prior discussed here.
Given a positive number $ C_1 $, for all $ \theta \in B_1 (C_1):= \{ \theta \in \mathbb{R}^d : \|\theta\|_1 \leq C_1 \} $, we consider the following prior,
\begin{eqnarray}
\label{eq_priordsitrbution}
\pi_\theta (\theta) 
\propto 
\prod_{i=1}^{d} 
	(\tau^2 + \theta_{i}^2)^{-2}
,
\end{eqnarray}
where $ \tau>0 $ is a tuning parameter.  This prior has been utilized in a variety of sparse contexts as in \cite{dalalyan2012mirror, dalalyan2012sparse,mai2023high}. Here, $C_1$ serves as a regularization parameter, typically assumed to be very large. Consequently, the distribution of $\pi_\theta$ closely resembles that of $S\tau \sqrt{2}$, where $S$ is a random vector with i.i.d. components drawn from a Student's t-distribution with $ 3$ degrees of freedom. By selecting an exceptionally small value for $\tau$, specifically less than $1/n$ as detailed later, the majority of elements in $\tau S$ concentrate near zero. However, due to the heavy-tailed nature, a small number of components in $\tau S$ deviate significantly from zero. This unique characteristic empowers the prior to promote sparsity in the parameter vector.
It is noted that the important of heavy-tailed priors in addressing sparsity has been studied before, as in \citep{seeger2008bayesian,johnstone2004needles,rivoirard2006nonlinear,abramovich2007optimality,carvalho2010horseshoe,castillo2012needles,castillo2015bayesian,castillo2018empirical,ray2022variational}.

We put an inverse-Gamma prior distribution on the  variance $ \sigma^2 $, $ \sigma^2 \sim IG (a,b) $. Here, we consider a fixed $ a>0 $ and $ b $ is a small tuning parameter.

\section{Results for fractional posterior}
\label{sc_concentration}
In this section we provide unified results regarding consistency as well as concentration rates of the fractional posterior in high-dimensional linear regression under suitable assumptions on the design matrix. Results for the regular posterior are presented in Section \ref{sc_regular_poster_n_horseshoe}.

\subsection{Assumptions}

\begin{assume}
	\label{asume_finite_2ndmoment}
	  Let assume that $ \mathbb{E}\left\Vert X_1\right\Vert^2 < C_{\rm x}^2 < \infty $.
\end{assume}

Assumption \ref{asume_finite_2ndmoment} is our main assumption to obtain our results, this assumption has been also used in \cite{jiang2007bayesian} for regression problem with random design.

\begin{assume}
	\label{asume_sigma_lowerbound}
Let's suppose that $ 0 < ( 2 \sigma_{\text{min}})^{-1} < \sigma_0^2 < ( 2 \sigma_{\text{max}})^{-1} < \infty $.
\end{assume}

Assumption \ref{asume_sigma_lowerbound} shares similarities with an assumption presented in  \citep{fang2024high}.

\begin{assume}
	\label{asume_minimaleigenvalue_estimebound}
	Assume that the minimal eigenvalue of the matrix $ G:= \mathbb{E}(XX^\top) $, denoted by $ \lambda_{min}^G $, is strictly positive.
\end{assume}
The above assumption has been previously employed in \cite{abramovich2018high}. For instance, it holds true when all columns $ X_{\cdot j} $ are linearly independent.

\subsection{Consistency results}
\label{sc_consistency}

Initially, a result in expectation is outlined, typically known as the consistency outcome of the fractional posterior.

\begin{theorem}
	\label{theorem_result_dis_expectation}
	For any $\alpha\in(0,1)$, 
assume that Assumption \ref{asume_finite_2ndmoment}, \ref{asume_sigma_lowerbound} hold. Put $ \tau = \frac1{C_{\rm x} \sqrt{nd}}\, $, for all $ \theta_0 $ such that $ \| \theta_0 \|_1 \leq C_1 - 2d\tau $, then
	\begin{equation*}
	\mathbb{E} \left[ \int D_{\alpha}(P_{\theta,\sigma} , P_{\theta_0,\sigma_0}) 
	\Pi_{n,\alpha}({\rm d}\theta , {\rm d} \sigma ) \right]
	\leq 
	\frac{1+\alpha}{1-\alpha}\varepsilon_n
	,
	\end{equation*}
	where
$
\varepsilon_n
=
c s^* \log 	\left(\frac{  C_{\rm x} C_1  \sqrt{nd} }{ s^*}\right) / n
,
$
for some positive constant $ c $ depending only on $ \sigma_{min},\sigma_{max},a $.
\end{theorem}

Using results from \cite{van2014renyi} regarding the connections between the Hellinger distance or the total variation and the $\alpha$-R\'enyi divergence, we derive the following results. Put
$ \kappa_\alpha := 2(\alpha+1) / (1-\alpha) $ when $ \alpha \in [0.5,1) $ and $ \kappa_\alpha := 2(\alpha+1) / \alpha $ when $ \alpha \in (0, 0.5) $.

\begin{cor}
	\label{cor_Hellinger_consistency}
	As a special case, Theorem \ref{theorem_result_dis_expectation} leads to a concentration result in terms of the classical Hellinger distance
	\begin{equation*}
	\mathbb{E}
	\left[
	\int H^2(P_{\theta,\sigma} ,P_{\theta_0,\sigma_0} ) \Pi_{n,\alpha}({\rm d}\theta , {\rm d} \sigma )  
	\right]
	\leq 
	\kappa_\alpha 
	\varepsilon_n 
	,
	\end{equation*}
	\begin{equation*}
	\mathbb{E}
	\left[	\int  d^{2}_{TV}(P_{\theta,\sigma} ,P_{\theta_0,\sigma_0} ) \Pi_{n,\alpha}({\rm d}\theta , {\rm d} \sigma ) 
	\right] 
	\leq 
	\frac{2(\alpha+1)}{(1-\alpha)\alpha}
	\varepsilon_n 
	,
	\end{equation*}
	with $d_{TV}$ being the total variation distance.
\end{cor}

As the squared Hellinger metric, total variation distance and $\alpha$-Rényi divergence entail inherent ambiguity, no claim is posited concerning the proximity between $ \theta $ and $ \theta_0 $ within a Euclidean-type distance framework. Deviating markedly from the emphasis of \cite{jiang2007bayesian,liang2013bayesian}, which solely relied on results derived from the Hellinger metric, our aim is to establish consistency outcomes for $ \theta $ utilizing a more clearly delineated metric.

\begin{cor}
\label{cor_consistency_in_ell2}
	 
	There exist some universal constant $ K>0 $ such that
	\begin{align*}
	\mathbb{E}
	\left[
	\int
	K (\mathbb{E}_{\rm x} [X_1^\top(\theta-\theta_0 )]^2
	+
	(\sigma- \sigma_0)^2 
	) 
	\Pi_{n,\alpha}({\rm d}\theta , {\rm d} \sigma )  
	\right]
	\leq 
	\kappa_\alpha 
	\varepsilon_n 
	,
	\end{align*}
and subsequently,
	\begin{align*}
	\mathbb{E}	\left[	\int
 \mathbb{E}_{\rm x} [X_1^\top(\theta-\theta_0 )]^2	
\Pi_{n,\alpha}({\rm d}\theta , {\rm d} \sigma )  
	\right]
	\leq 
	 \frac{\kappa_\alpha}{K}
	\varepsilon_n 
	,
\\
	\mathbb{E}	\left[	\int
	(\sigma- \sigma_0)^2 
\Pi_{n,\alpha}({\rm d}\theta , {\rm d} \sigma )  
	\right]
	\leq 
	 \frac{\kappa_\alpha}{K}
	\varepsilon_n 
	,
\end{align*}
and under additional Assumption \ref{asume_minimaleigenvalue_estimebound} holds, we obtain that
	\begin{align*}
\mathbb{E}
\left[
\int
\| \theta-\theta_0 \|_2^2	
\Pi_{n,\alpha}({\rm d}\theta , {\rm d} \sigma )  
\right]
\leq 
 \frac{\kappa_\alpha}{\lambda_{min}^G K}
\varepsilon_n 
.
\end{align*}
\end{cor}

These results are novel up to the best of our knowledge.

\subsection{Concentration rates}
\label{sc_rslton_distribu}

We now present our principal finding regarding the concentration of the fractional posterior concerning the $\alpha$-Rényi divergence of the densities.

  \begin{theorem}
\label{theorem_main}
	For any $\alpha\in(0,1)$, 
assuming that Assumption \ref{asume_finite_2ndmoment}, \ref{asume_sigma_lowerbound} hold. Put $ \tau = \frac1{C_{\rm x} \sqrt{nd}}\, $.
For any $ (\varepsilon,\eta)\in(0,1)^2 $,
 and for all $ \theta_0 $ that $ \|\theta_0\|_1 \leq C_1 - 2d\tau  $,
     we have that
    \begin{equation*}
    \mathbb{P}\left[
 \int D_{\alpha}(P_{\theta,\sigma} ,P_{\theta_0,\sigma_0} ) \Pi_{n,\alpha}({\rm d}\theta , {\rm d} \sigma )  
 \leq  
 \frac{(\alpha+1) \varepsilon_n + \alpha \sqrt{\frac{\varepsilon_n}{n\eta}} 
 	+ 
 	\frac{
 	\log\left( 1/\varepsilon \right)}{n}}{1-\alpha}\right]
 \geq 
 1-\varepsilon-\eta.
 \end{equation*}
	where
$
\varepsilon_n
=
c  s^* \log 
	\left(\frac{  C_{\rm x} C_1  \sqrt{nd} }{ s^*}\right) /n
,
$
for some positive constant $ c $ depending only on $  \sigma_{min},\sigma_{max},a $.
  \end{theorem}

The argument in the proof of Theorem \ref{theorem_main} adheres to a general approach for the fractional posterior as outlined in \cite{bhattacharya2016bayesian} and \cite{alquier2020concentration}. Notably, the concentration of the fractional posterior is achieved solely through the properties of prior concentration, as discussed immediately after Theorem 2.4 in \cite{alquier2020concentration}. In contrast, the concentration theory for the regular posterior requires more stringent conditions, as highlighted in \cite{jeong2021posterior} and also discussed in \cite{jiang2007bayesian}. This observation has been emphasized in various studies, including \cite{bhattacharya2016bayesian, alquier2020concentration, cherief2018consistency, chakraborty2020bayesian, l2023semiparametric}.

Following \cite{alquier2020concentration}, taking $ \eta = (n \varepsilon_n)^{-1} $ and $\varepsilon = \exp(-n \varepsilon_n)$, one obtain a concentration result by noting that $ 1- (n \varepsilon_n)^{-1} - \exp(-n\varepsilon_n) \geq 
1- 2(n \varepsilon_n)^{-1} $. Therefore, the sequence $\varepsilon_n$ gives a concentration rate for the fractional posterior in \eqref{eq_fractional_posterior}, stated in the following corollary. 
\begin{cor}
	\label{cor_concentration}
	Under the same assumptions as in Theorem~\ref{theorem_main},
	\begin{align*}
	\mathbb{P}\left[
	\int D_{\alpha}(P_{\theta,\sigma} ,P_{\theta_0,\sigma_0} ) \Pi_{n,\alpha}({\rm d}\theta , {\rm d} \sigma )  
	\leq 
	\frac{2(\alpha+1)}{1-\alpha} \varepsilon_n\right] 
	\geq 
	1-\frac{2}{n\varepsilon_n}
	.
	\end{align*}
\end{cor}

\begin{rmk}
	The main assumption regarding the distribution of $X_1$, in both Theorem \ref{theorem_result_dis_expectation} and Theorem \ref{theorem_main}, is that its second moment is bounded and this does not effect our concentration rates. For example, when $X_1$ follows a uniform distribution on the unit sphere, then $ \mathbb{E}\left\Vert X_1\right\Vert^2 \leq 4$. When $X_1 \sim \mathcal{N}(0, \vartheta^2 I_d)$, the values are $ \mathbb{E}\left\Vert X_1\right\Vert^2 = 4 \vartheta^2 d$.
\end{rmk}

\begin{rmk}
We remind that all technical proofs are provided in Appendix \ref{sc_proofs}. Our results demonstrate that the concentration rates of the fractional posterior adapt to the unknown sparsity level \( s^* \). While prior studies, such as \cite{jiang2007bayesian, liang2013bayesian, jeong2021posterior}, focus on contraction rates based on the Hellinger distance, our work distinguishes itself by offering comprehensive concentration results for the fractional posterior using the \(\alpha\)-Rényi divergence. Specifically, we derive the following notable corollary on the Hellinger distance and the total variation distance, leveraging results from \cite{van2014renyi}.
\end{rmk}

\begin{cor}
	\label{cor_concentration_Hellinger}
	As a special case, Theorem \ref{theorem_main} leads to a concentration result in terms of the classical Hellinger distance
	\begin{equation*}
	\mathbb{P}\left[
	\int H^2(P_{\theta,\sigma} ,P_{\theta_0,\sigma_0} ) \Pi_{n,\alpha}({\rm d}\theta , {\rm d} \sigma )  
	\leq 
	\kappa_\alpha 
	\varepsilon_n \right]
	\geq 
	1-\frac{2}{n\varepsilon_n}
	.
	\end{equation*}
	And for $ \alpha \in (0,1) $, 
	\begin{equation*}
	\mathbb{P}\left[
	\int  d^{2}_{TV}(P_{\theta,\sigma} ,P_{\theta_0,\sigma_0} ) \Pi_{n,\alpha}({\rm d}\theta , {\rm d} \sigma )  
	\leq 
	\frac{4(\alpha+1)}{(1-\alpha)\alpha}
	\varepsilon_n \right]
	\geq 
	1-\frac{2}{n\varepsilon_n}
	,
	\end{equation*}
 with $d_{TV}$ being the total variation distance.
\end{cor}

Corollary \ref{cor_concentration_Hellinger} demonstrates that the fractional posterior distribution of \( \theta \) and \( \sigma \) concentrates around its true value at a specified rate relative to the squared Hellinger metric and to the total variation distance. However, the implications of our findings differ from the most recent work by \cite{jeong2021posterior}: we consider random design and unknown error variance with i.i.d. observations, whereas \cite{jeong2021posterior} dealt with independent observations for fixed design and known error variance.

Given the ambiguity inherent in the squared Hellinger metric, the total variation, and the \(\alpha\)-Rényi divergence, we refrain from making assertions about the proximity of \( \theta \) and \( \theta_0 \) in the context of a Euclidean-type distance. Our objective now is to furnish concentration results for parameters utilizing a more explicitly defined metric.

\begin{cor}
	\label{cor_concen_rate_in_ell2}
	There exist some universal constant $ K>0 $ such that
	\begin{align*}
	\mathbb{P}	\left[
	\int
	K (\mathbb{E}_{\rm x} [X_1^\top(\theta-\theta_0 )]^2
	+
	(\sigma- \sigma_0)^2 
	) 
	\Pi_{n,\alpha}({\rm d}\theta , {\rm d} \sigma )  
	\leq 
\kappa_\alpha 
\varepsilon_n \right]
\geq 
1-\frac{2}{n\varepsilon_n}
	,
	\end{align*}
	and subsequently,
	\begin{align*}
	\mathbb{P}
	\left[	\int
	 \mathbb{E}_{\rm x} [X_1^\top(\theta-\theta_0 )]^2	
	\Pi_{n,\alpha}({\rm d}\theta , {\rm d} \sigma )  
	\leq 
 \frac{\kappa_\alpha}{K}
\varepsilon_n \right]
\geq 
1-\frac{2}{n\varepsilon_n}
	,
	\end{align*}
	\begin{align*}
\mathbb{P}
\left[	\int
	(\sigma- \sigma_0)^2 
\Pi_{n,\alpha}({\rm d}\theta , {\rm d} \sigma )  
\leq 
\frac{\kappa_\alpha}{K}
\varepsilon_n \right]
\geq 
1-\frac{2}{n\varepsilon_n}
,
\end{align*}
	and under additional Assumption \ref{asume_minimaleigenvalue_estimebound} holds, we obtain that
	\begin{align*}
	\mathbb{P}
	\left[	\int
	\| \theta-\theta_0 \|_2^2	
	\Pi_{n,\alpha}({\rm d}\theta , {\rm d} \sigma )  
	\leq 
 \frac{\kappa_\alpha}{\lambda_{min}^G K}
\varepsilon_n \right]
\geq 
1-\frac{2}{n\varepsilon_n}
	.
	\end{align*}
\end{cor}

\begin{rmk}
The prediction error, $ \| X_1^\top(\theta-\theta_0) \|_2^2 $, and the estimation error, $ \| \theta-\theta_0 \|_2^2 $, of order $ s^*\log(d/s^*)/n $, are minimax-optimal and are similar to those found in the frequentist analysis by \cite{bellec2018slope} and \cite{derumigny2018improved}. However, the results in \cite{bellec2018slope} and \cite{derumigny2018improved} require that the random design satisfies a type of Restricted Eigenvalue condition, which our approach does not necessitate.
\end{rmk}

\subsection{Result in the misspecified case}
\label{sc_mispecifiedcase}
In this section, we show that our results in Section \ref{sc_consistency} can be extended to the misspecified setting.
Assume that the true data generating distribution is parametrized by $ \theta_0 \notin \Theta $ and define $P_{\theta_0,\sigma_0} $ as the true distribution. 
Put
\begin{align*}
\theta^* 
:=
\arg\min_{\theta \in \Theta} 
\mathcal{K}(P_{\theta_0,\sigma_0} ,P_{\theta,\sigma_0} ).
\end{align*}
In order not to change all
the notation, we consider an extended parameter set as $ \{\theta_0 \} \cup \Theta$ and obtain the following result.

\begin{theorem}
	\label{theorem_misspecified}
	For any $\alpha\in(0,1)$, let assume that Assumption \ref{asume_finite_2ndmoment} hold, $ \| \theta_0 \|_1 \leq C_1 - 2d\tau $, and with $ \tau = (C_{\rm x} \sqrt{nd})^{-1} $. Then, 
	\begin{equation*}
	\mathbb{E} \left[ \int D_{\alpha}(P_{\theta,\sigma} ,P_{\theta_0,\sigma_0} ) \Pi_{n,\alpha}({\rm d}\theta , {\rm d} \sigma )\right]
	\leq 
	\frac{\alpha}{1-\alpha} \min_{\theta\in\Theta} \mathcal{K}(P_{\theta_0,\sigma_0} ,P_{\theta,\sigma_0} )
	+ 
	c\frac{1+\alpha}{1-\alpha} \varepsilon_n
	,
	\end{equation*}
	where
$
\varepsilon_n
=
 s^* \log 
	\left(\frac{  C_{\rm x} C_1  \sqrt{nd} }{ s^*}\right)/n
,
$
for some constant $ c>0 $ depending only on $  \sigma_{min},a $.
\end{theorem}

In the case of a well-specified model, i.e., when $\theta_0 = \arg\min_{\theta\in\Theta} \mathcal{K}(P_{\theta_0,\sigma_0} ,P_{\theta,\sigma} )$, one recovers Theorem~\ref{theorem_result_dis_expectation}, otherwise, this result presents an oracle inequality. Although it may not constitute a sharp oracle inequality due to the differing risk measures on both sides, this observation remains valuable, particularly when $\mathcal{K}(P_{\theta_0,\sigma_0} ,P_{\theta^*,\sigma_0})$ is minimal. Nonetheless, under additional assumptions, we can further derive an oracle inequality result with $\ell_2$ distance on both sides. The result is as follows.

\begin{cor}
	\label{cor_misspecified}
Assume that Theorem \ref{theorem_misspecified} holds and additional Assumption \ref{asume_sigma_lowerbound} is satisfied. Then, for some positive constant $ K_\alpha $ depending only on $ \alpha $,
 	\begin{equation*}
\mathbb{E} \left[ \int
\mathbb{E}_{\rm x} [X_1^\top(\theta-\theta_0 )]^2
\Pi_{n,\alpha}({\rm d}\theta , {\rm d} \sigma )\right]
\leq 
K_\alpha 
\left\{
\min_{\theta\in\Theta} \sigma_{min} \mathbb{E}_{\rm x} [X_1^\top(\theta-\theta_0 )]^2
+  \varepsilon_n
\right\}
,
\end{equation*}
and under Assumption \ref{asume_minimaleigenvalue_estimebound},
 	\begin{equation*}
\mathbb{E} \left[ \int
\lambda_{min}^G \| \theta - \theta_0 \|^2
\Pi_{n,\alpha}({\rm d}\theta , {\rm d} \sigma )\right]
\leq 
K_\alpha 
\left\{
\min_{\theta\in\Theta} \sigma_{min} 
C_{\rm x}^2 \| \theta - \theta_0 \|^2
+  \varepsilon_n
\right\}
.
\end{equation*}
\end{cor}

\begin{rmk}
While oracle inequalities for penalized procedures have been established by \cite{bellec2018slope}, our findings in Corollary \ref{cor_misspecified} and Corollary \ref{cor_regular_poster_in_ell2} (pertaining to the regular posterior) represent  novel contributions for Bayesian approaches in the context of high-dimensional sparse linear regression with random design, to the best of our knowledge. It is important to note that our oracle inequalities are not sharp, as the leading constants on the right-hand side are not equal to 1, unlike those in \cite{bellec2018slope}. Additionally, non-sharp oracle inequalities for Bayesian methods have been derived by \cite{castillo2015bayesian}, but these pertain to the fixed design setting.
\end{rmk}

\section{Result for regular posterior and for other priors}
\label{sc_regular_poster_n_horseshoe}
\subsection{Result for regular posterior}
\label{sc_regular_posterior}
In this section, we provide results concerning the regular posterior distribution, for which $ \alpha =1 $ in \eqref{eq_fractional_posterior}.

For any $\alpha \in (0,1)$, one has that
\begin{align}
\Pi_{n}(\theta,\sigma )
& 
\propto 
\sigma^{-n} \, e^{- \sum_{i=1}^{n} \frac{(Y_i - X_i^\top \theta)}{2\sigma^2} } 
 \pi_\theta ({\rm d}\theta) \pi_\sigma({\rm d}\sigma)
 \nonumber
\\
& \propto 
\sigma^{-n \alpha} \, e^{- \alpha  
\sum_{i=1}^{n} \frac{(Y_i - X_i^\top \theta)}{\alpha 2\sigma^2} }  
 \pi_\theta ({\rm d}\theta)
 \, \sigma^{-n (1-\alpha)} \,  \pi_\sigma({\rm d}\sigma)
 \nonumber
\\
& \propto 
\sigma_*^{-n\alpha} \, e^{- \alpha \sum_{i=1}^{n} \frac{(Y_i - X_i^\top \theta)}{ 2\sigma_*^2} } 
\, \pi_\theta ({\rm d}\theta)
 \pi_{\sigma_*}({\rm d}\sigma_*)
 \nonumber
\\
& \propto 
\Pi_{n,\alpha}( \theta, \sigma_* )
\label{eq_connection_fracpost_regu_post}
,
\end{align} 
where $\sigma_* = \alpha \sigma$ and $ \pi_{\sigma_*}(\cdot)  $ is again an $ IG (n(1-\alpha) + a, \alpha b) $. Given that the first and last expressions in \eqref{eq_connection_fracpost_regu_post} both represent probability densities, it follows that \( \Pi_{n}(\theta,\sigma ) = \Pi_{n,\alpha}( \theta, \sigma_* ) \). Similar relationships have been acknowledged in different contexts, as in \cite{chakraborty2020bayesian,martin2020empirical}.

This implies that the regular posterior distribution of $(\theta, \sigma)$ can be interpreted as the $\alpha$-fractional posterior distribution of $(\theta, \sigma_*)$, where the prior distribution of $\sigma_*$ is contingent on both $n$ and $\alpha$. Following an argument similar to Theorem \ref{theorem_main}, we only need to show the prior concentration of $( \theta, \sigma_*)$ around the truth to obtain posterior consistency of $ \Pi_{n,\alpha}( \theta, \sigma_* ) $, and hence equivalently of $\Pi_{n}( \theta, \sigma ) $. The only place that needs additional care is showing the prior concentration of $\sigma_*$ with an $n$-dependent prior. The prior for $ \theta $ remains the same as in Section \ref{sc_concentration}.

\begin{theorem}
\label{thm_regular_posterior}
Assume that Assumption \ref{asume_finite_2ndmoment}, \ref{asume_sigma_lowerbound} hold. Put $ \tau = \frac1{C_{\rm x} \sqrt{nd}}\, $
and for all $ \theta_0 $ that $ \|\theta_0\|_1 \leq C_1 - 2d\tau  $,
we have that
	\begin{equation*}
\mathbb{E} \left[ \int D_{\alpha}(P_{\theta,\sigma} , P_{\theta_0,\sigma_0}) 
\Pi_{n}({\rm d}\theta , {\rm d} \sigma ) \right]
\leq 
\frac{1+\alpha}{1-\alpha}\varepsilon_n
,
\end{equation*}
and
	\begin{align*}
	\mathbb{P}\left[
	\int D_{\alpha}(P_{\theta,\sigma} ,P_{\theta_0,\sigma_0} ) \Pi_{n } 
	({\rm d}\theta , {\rm d} \sigma )  
	\leq 
	\frac{2(\alpha+1)}{1-\alpha} \varepsilon_n\right] 
	\geq 
	1-\frac{2}{n\varepsilon_n}
	,
	\end{align*}
	where
$
	\varepsilon_n
	=
	c s^* \log 	\left(\frac{ C_{\rm x} C_1 \sqrt{nd} }{ s^*}\right) / n
	,
$
	for some constant $ c>0 $ depending only on $  \sigma_{min},\sigma_{max}, a $.
\end{theorem}

Theorem \ref{thm_regular_posterior} establishes the consistency and concentration rates of the regular posterior relative to the \(\alpha\)-R\'enyi divergence. As \(n < d\), the concentration rate is of the order $ s^* \log 	\left(  d/s^* \right) / n $. Additionally, the upper bound constants \(C_1\) and \(C_{\rm x}\) can be relaxed to grow polynomially in terms of \(n\) and \(d\) without affecting the rate.

In line with Corollary \ref{cor_concentration_Hellinger} for the fractional posterior, we immediately obtain the following results from Theorem \ref{thm_regular_posterior}. The proof follows a similar structure and is thus omitted.

\begin{cor}
	\label{cor_concen_Hellinger_regular_posterior}
	As a special case, Theorem \ref{theorem_main} leads to a concentration result in terms of the classical Hellinger distance
	\begin{equation*}
	\mathbb{P}\left[
	\int H^2(P_{\theta,\sigma} ,P_{\theta_0,\sigma_0} ) 
	\Pi_{n}({\rm d}\theta , {\rm d} \sigma )  
	\leq 
	\kappa_\alpha 
	\varepsilon_n \right]
	\geq 
	1-\frac{2}{n\varepsilon_n}
	,
	\end{equation*}
	\begin{equation*}
	\mathbb{P}\left[
	\int  d^{2}_{TV}(P_{\theta,\sigma} ,P_{\theta_0,\sigma_0} ) 
	\Pi_{n}({\rm d}\theta , {\rm d} \sigma )  
	\leq 
	\frac{4(\alpha+1)}{(1-\alpha)\alpha}
	\varepsilon_n \right]
	\geq 
	1-\frac{2}{n\varepsilon_n}
	,
	\end{equation*}
	with $d_{TV}$ being the total variation distance.
\end{cor}

Our results imply that the concentration rates of the regular posterior are adaptive to the unknown sparsity level $ s^* $ with random design and unknown error variance.

In accordance with Corollary \ref{cor_concen_rate_in_ell2} for the fractional posterior, we immediately derive the following results from Theorem \ref{thm_regular_posterior}. The proof follows a similar structure and is thus omitted.

\begin{cor}
	\label{cor_regular_poster_in_ell2}
	There exist some universal constant $ K>0 $ such that
	\begin{align*}
	\mathbb{P}	\left[
	\int
	K (\mathbb{E}_{\rm x} [X_1^\top(\theta-\theta_0 )]^2
	+
	(\sigma- \sigma_0)^2 
	) 
	\Pi_{n}({\rm d}\theta , {\rm d} \sigma )  
	\leq 
	\kappa_\alpha 
	\varepsilon_n \right]
	\geq 
	1-\frac{2}{n\varepsilon_n}
	,
	\end{align*}
	and subsequently,
	\begin{align*}
	\mathbb{P}
	\left[	\int
	K \mathbb{E}_{\rm x} [X_1^\top(\theta-\theta_0 )]^2	
	\,
	\Pi_{n}({\rm d}\theta , {\rm d} \sigma )  
	\leq 
	\kappa_\alpha 
	\varepsilon_n \right]
	\geq 
	1-\frac{2}{n\varepsilon_n}
	,
	\end{align*}
	and under additional Assumption \ref{asume_minimaleigenvalue_estimebound} holds, we obtain that
	\begin{align*}
	\mathbb{P}
	\left[	\int
	\| \theta-\theta_0 \|_2^2	
	\,
	\Pi_{n}({\rm d}\theta , {\rm d} \sigma )  
	\leq 
 \frac{\kappa_\alpha}{\lambda_{min}^G K}
	\varepsilon_n \right]
	\geq 
	1-\frac{2}{n\varepsilon_n}
	.
	\end{align*}
\end{cor}

\begin{rmk}
	It is noteworthy to note that the assumption outlined in Assumption \ref{asume_minimaleigenvalue_estimebound} on the minimal eigenvalue of the gram matrix $ G $, can be substituted with broader conditions. Specifically, we can leverage the compatibility numbers, defined as
	$$ \phi_1 (s) 
	:= 
	\inf_{\theta:0 < \|\theta\|_0 \leq s} 
	\frac{ \mathbb{E}_{\rm x} \| X^\top\theta \|_2^2 \|\theta\|_0 }{ \|\theta \|_1^2 }
	,
	\quad
	\phi_2 (s) 
	:= 
	\inf_{\theta:0 < \|\theta\|_0 \leq s} 
	\frac{ \mathbb{E}_{\rm x} \| X^\top\theta \|_2^2 }{  \|\theta \|_2^2 }
	$$
	to derive results in $\ell_1$ and $\ell_2$ distances. Analogous to Corollary \ref{cor_regular_poster_in_ell2}, we find that
		\begin{align*}
	\mathbb{P}
	\left[	\int
	\| \theta-\theta_0 \|_2^2	
	\,
	\Pi_{n}({\rm d}\theta , {\rm d} \sigma )  
	\leq 
	\frac{\kappa_\alpha}{	\phi_2 (s^*)  K}
	\varepsilon_n \right]
	\geq 
	1-\frac{2}{n\varepsilon_n}
	.
	\end{align*}
		\begin{align*}
	\mathbb{P}
	\left[	\int
	\| \theta-\theta_0 \|_1^2	
	\,
	\Pi_{n}({\rm d}\theta , {\rm d} \sigma )  
	\leq 
	\frac{\kappa_\alpha}{	\phi_1 (s^*)  K}
s^*	\varepsilon_n \right]
	\geq 
	1-\frac{2}{n\varepsilon_n}
	.
	\end{align*}
The notion of compatibility numbers is frequently utilized in Bayesian high-dimensional literature, as illustrated in studies such as \cite{castillo2015bayesian, martin2017empirical, ray2020spike, belitser2020empirical, jeong2021posterior}. Originally emerging from high-dimensional frequentist literature, compatibility numbers were first employed in contexts like those discussed in \cite{buhlmann2011statistics}. The relationship between compatibility numbers and the Restricted Eigenvalue condition is explored in \cite{bellec2018slope}.
\end{rmk}

\subsubsection{Result in the misspecified case}
We now study the property of the regular posterior in a model misspecified case.
Assume that the true data generating distribution $P_{\theta_0,\sigma_0} $ is now parametrized by $ \theta_0 \notin \Theta $.
Put
\begin{align*}
\theta^* 
:=
\arg\min_{\theta \in \Theta} 
\mathcal{K}(P_{\theta_0,\sigma_0} ,P_{\theta,\sigma_0} ),
\end{align*}
For the sake of simplicity, we define an extended parameter set $ \{\theta_0 \} \cup \Theta$, we obtain the following result.

\begin{theorem}
	\label{thm_misspecified_regularposte}
Assume that Assumption \ref{asume_finite_2ndmoment} hold, $ \| \theta_0 \|_1 \leq C_1 - 2d\tau $, and with $ \tau = (C_{\rm x} \sqrt{nd})^{-1} $. Then, 
	\begin{equation*}
	\mathbb{E} \left[ \int D_{\alpha}(P_{\theta,\sigma} ,P_{\theta_0,\sigma_0} ) \Pi_{n}({\rm d}\theta , {\rm d} \sigma )\right]
	\leq 
	\frac{\alpha}{1-\alpha} \min_{\theta\in\Theta} \mathcal{K}(P_{\theta_0,\sigma_0} ,P_{\theta,\sigma_0} )
	+ \frac{1+\alpha}{1-\alpha} \varepsilon_n
	,
	\end{equation*}
	where
	$
	\varepsilon_n
	=
	c
	s^* \log 
	\left(\frac{  C_{\rm x} C_1  \sqrt{nd} }{ s^*}\right)/n
	,
	$
	for some constant $ c>0 $ depending only on $  \sigma_{min},\sigma_{max},a $.
\end{theorem}

In the case of a well-specified model, i.e., when $\theta_0 = \arg\min_{\theta\in\Theta} \mathcal{K}(P_{\theta_0,\sigma_0} ,P_{\theta,\sigma_0} )$, one recovers Theorem~\ref{thm_regular_posterior}, otherwise, this result presents an oracle inequality. Although it it not a sharp oracle inequality due to the differing risk measures on both sides, this observation remains valuable, particularly when $\mathcal{K}(P_{\theta_0,\sigma_0} ,P_{\theta^*,\sigma_0})$ is minimal. 

To the best of our knowledge, the result in Theorem \ref{thm_misspecified_regularposte} regarding the regular posterior is completely novel in the context of sparse high-dimensional linear regression with random design and unknown error variance.

Nevertheless, under additional assumptions, we can derive an oracle inequality result with \(\ell_2\) distance on both sides. In line with Corollary \ref{cor_misspecified} for the fractional posterior, we immediately obtain the following results from Theorem \ref{thm_misspecified_regularposte}. The proof follows a similar structure and is thus omitted.

\begin{cor}
	\label{cor_misspecified_regularposte}
	Assume that Theorem \ref{thm_misspecified_regularposte} holds and additional Assumption \ref{asume_sigma_lowerbound} is satisfied. Then, 
	\begin{equation*}
	\mathbb{E} \left[ \int
	\mathbb{E}_{\rm x} [X_1^\top(\theta-\theta_0 )]^2
	\Pi_{n}({\rm d}\theta , {\rm d} \sigma )\right]
	\leq 
	K_\alpha 
	\left\{
	\min_{\theta\in\Theta} \sigma_{min} \mathbb{E}_{\rm x} [X_1^\top(\theta-\theta_0 )]^2
	+  \varepsilon_n
	\right\}
	.
	\end{equation*}
\end{cor}

\subsection{Result with spike and slab prior}
\label{sc_spike_n_slab_prior}

The use of the scaled Student's t-distribution, as discussed in the previous section, has shown its advantages. Nonetheless, it is not suitable for variable selection. In situations where sparsity is important, the spike and slab prior, as proposed by \cite{mitchell1988bayesian,george1993variable}, may be preferable for variable selection purposes.

\begin{equation}
\label{eq_spike_slab_prior}
\pi_{SnS}(\theta) 
= 
\prod_{i=1}^{d} \left[p \phi (\theta_i;0,v_1)
+ (1-p) \phi(\theta_i;0,v_0) \right] 
\end{equation}
with $ (p,v_0,v_1)\in[0,1]\times (\mathbb{R}^+)^2$, and $v_0 \ll v_1$. Here, $ \phi (\cdot ; m, s) $ is the Gaussian density with mean $ m $ and variance $ s $. As \( v_0 \) tends to $ 0 $, it converges to a more traditional prior that places a point mass at zero for each component. However, this leads to a (fractional) posterior consisting of a mixture of \(2^d\) components, combining Dirac masses and continuous distributions, thereby making the computation more complex.

Utilizing the spike and slab prior as described in equation \eqref{eq_spike_slab_prior} can produce results similar to those obtained using the scaled Student prior outlined in equation \eqref{eq_priordsitrbution}. Theorem \ref{theorem_spike_slab}, presented below, provides outcomes comparable to those derived from the main result, Theorem \ref{theorem_main}.

  \begin{theorem}
	\label{theorem_spike_slab}
	Let assume that Assumption \ref{asume_finite_2ndmoment} hold and $\|\theta_0\|_2 \leq 1 $. Using spike and slab prior in \eqref{eq_spike_slab_prior} with $ p = 1-e^{-1/d} $, $ v_0 \leq 1/(2nd\log(2)) $. 
Then, 
	\begin{align*}
\mathbb{P}\left[
\int D_{\alpha}(P_{\theta,\sigma} ,P_{\theta_0,\sigma_0} ) \Pi_{n,\alpha}({\rm d}\theta , {\rm d} \sigma )  
\leq 
c_{v_1}\frac{2(\alpha+1)}{1-\alpha} \varepsilon_n\right] 
\geq 
1-\frac{2}{n\varepsilon_n}
,
\end{align*}
where $ c_{v_1} >0 $ is a constant depending only on $ v_1 $, and
$
\varepsilon_n
=
s^* \log 	\left(\frac{ C_{\rm x} \sqrt{nd} }{ s^*}\right) / n
.
$
\end{theorem}

It is worth mentioning that in \cite{rovckova2014emvs}, the authors successfully reduced the computational burden associated with aspects of the posterior distribution by replacing the spike at zero with a Gaussian distribution with a small variance as in \eqref{eq_spike_slab_prior}. More specifically, they proposed an efficient EM algorithm for this purpose. Thus, our Theorem \ref{theorem_spike_slab} presents a theoretical perspective for the method proposed in \cite{rovckova2014emvs}. The proof is given in Appendix \ref{sc_proofs}.

\section{Conclusion}
\label{sc_discuconcluse}
In this paper, we have examined the consistency and concentration rates of both the regular posterior and the fractional posterior in high-dimensional sparse linear regression models with random design and unknown error variance. By leveraging techniques from the fractional posterior, we have obtained novel results for the regular posterior in the context of model misspecification. Our results apply to two different types of priors: continuous shrinkage priors and spike-and-slab priors. A potential gap in our work is the lack of a sharp oracle inequality for the model misspecified case, similar to those found in the frequentist literature; this remains an open question for future research.

\subsubsection*{Acknowledgments}
The author acknowledges support from  Norwegian Research Council grant no. 309960, through the Centre for Geophysical Forecasting at NTNU. 

\subsubsection*{Conflicts of interest/Competing interests}
The author declares no potential conflict of interests.

\appendix
\section{Proofs}
\label{sc_proofs}

The following useful result concerning the Kullback-Leibler divergence between two Gaussian distributions is crucial for our proofs,
\begin{align*}
\mathcal{K} (\mathcal{N}(\mu_1,\sigma_1^2),\mathcal{N}(\mu_2,\sigma_2^2)) 
=
\log\left(\frac{\sigma_2}{\sigma_1}\right) + \frac{\sigma_1^2 + (\mu_1-\mu_2)^2}{2\sigma_2^2} - \frac{1}{2}
.
\end{align*}

\subsection{Proof of Section \ref{sc_concentration}}
\subsubsection{Proofs for Subsection \ref{sc_consistency}}
\begin{proof}[\bf Proof of Theorem~\ref{theorem_result_dis_expectation}]
	
	We will utilize the findings outlined in Theorem 2.6 of  \cite{alquier2020concentration}. To proceed, it is necessary to confirm that the conditions specified therein are satisfied.

First, we have 
 \begin{align*}
 	\mathcal{K}(P_{\theta_0,\sigma_0} ,P_{\theta,\sigma} )
& = 
\frac{ \mathbb{E}_{\rm x} (X_1^\top(\theta-\theta_0 )^2) }{2\sigma^2} 
+
\log\left(\frac{\sigma}{\sigma_0}\right) 
+ 
\frac{\sigma_0^2}{2\sigma^2} - \frac{1}{2}
\\
& \leq
 \frac{\mathbb{E} \| X_1 \|^2 \| \theta-\theta_0\|_2^2  }{2\sigma^2}
 +
 \log\left(\frac{\sigma}{\sigma_0}\right) 
 + 
 \frac{\sigma_0^2}{2\sigma^2} - \frac{1}{2}
 .
 \end{align*}
When integrating with respect to $\rho_n = p_0\otimes \tilde{\pi}_\sigma $, as defined in equations \eqref{eq_specific_distribution} and \eqref{eq_prior_sigma_specific}, and employing Lemma \ref{lema_boundfor_ell2} along with Lemma \ref{lm_bound_for_kl_term}, we obtain the following:
	\begin{align*}
	\int\mathcal{K}(P_{\theta_0,\sigma_0} ,P_{\theta,\sigma} )\rho_n({\rm d}\theta , {\rm d} \sigma )
&	\leq
\frac{2C_{\rm x}^2 }{\sigma_{min}^{-1}}
	\int \Vert\theta-\theta_0\Vert^2_2
p_0 ({\rm d}\theta , {\rm d} \sigma )
+ K  
\frac{\log(\varrho) \eta^2}{	\sigma_{min}^{-2}} 
\\
&	\leq 
\frac{2C_{\rm x}^2 }{\sigma_{min}^{-1}} 
4d \tau^2 
+ K  
\frac{\log(\varrho) \eta^2}{ \sigma_{min}^{-2} }
	,
	\end{align*} 
	where we have used Assumption \ref{asume_finite_2ndmoment}.
From Lemma \ref{lema_boundof_KL} and Lemma \ref{lm_sigma_prior}, we have that 
\begin{align*}
\frac1n \mathcal{K}(\rho_n,\pi)
 =
\frac1n \mathcal{K}(p_0,\pi_{\theta} )
+ 
\frac1n \mathcal{K} (\tilde{\pi}_\sigma, \pi_\sigma)
\leq
\frac{	4 s^* \log \left(\frac{C_1 }{\tau s^*}\right)
	+	\log(2)
 }{n}
+ 
\frac{
	\log \left( \frac{e 2^{a+1}\Gamma(a) }{ \eta^a} \right)
}{n}
.
\end{align*}
Putting $ \tau = \frac1{\sqrt{nd} C_{\rm x}}\, $ and $ \eta^2 = 1/n $, it leads to
	\begin{align*}
\int\mathcal{K}(P_{\theta_0,\sigma_0} ,P_{\theta,\sigma} )\rho_n({\rm d}\theta , {\rm d} \sigma )
&	\leq
\frac{8 }{n\sigma_{min}^{-1}}
+ K  
\frac{\log(\varrho)}{n \sigma_{min}^{-2} }
,
\end{align*}
and
\begin{align*}
\frac1n \mathcal{K}(\rho_n,\pi)
\leq
\frac{	4 s^* \log \left(\frac{C_{\rm x} C_1 \sqrt{nd} }{ s^*}\right)
+	\log(2)
+ \log \left( e 2^{a+1}\Gamma(a) n^{a/2} \right)}{n}
.
\end{align*}
Therefore, we can now apply Theorem 2.6 in \cite{alquier2020concentration} with $\rho_n $ and
 with
\[
\varepsilon_n
=
c s^* \log 
\left(\frac{ C_{\rm x} C_1 \sqrt{nd} }{ s^*}\right) /n
,
\]
for some constant $ c>0 $ depending only on $  \sigma_{min},\sigma_{max},a $.
 This completes the proof.
	
\end{proof}

\begin{proof}[\bf Proof of Corollary \ref{cor_Hellinger_consistency}]
	From \cite{van2014renyi}, we have that 
	$$ 
	H^2(P_{\theta,\sigma} ,P_{\theta_0,\sigma_0} )
	\leq 
	D_{1/2}(P_{\theta,\sigma} ,P_{\theta_0,\sigma_0} )  
	\leq 
	D_{\alpha}(P_{\theta,\sigma} ,P_{\theta_0,\sigma_0} ) 
	,
	$$
	for $ \alpha \in [0.5,1) $. In addition, 	for $ \alpha \in (0, 0.5)  $, we also have that 
	$$ 
	D_{1/2}(P_{\theta,\sigma} ,P_{\theta_0,\sigma_0} )  
	\leq 
	\frac{(1-\alpha)1/2}{\alpha (1-1/2)} D_{\alpha}(P_{\theta,\sigma} ,P_{\theta_0,\sigma_0} ) 
	=
	\frac{(1-\alpha)}{\alpha} D_{\alpha}(P_{\theta,\sigma} ,P_{\theta_0,\sigma_0} ) 
	.
	$$ 
	
	Thus, for
 \begin{equation*}
	\kappa_\alpha 
	=
	\begin{cases}
	\frac{2(\alpha+1)}{1-\alpha}, \alpha \in [0.5,1) ,
	\\
	\frac{2(\alpha+1)}{\alpha}, \alpha \in (0, 0.5)
	.
	\end{cases}
	\end{equation*} 
	and from  Theorem \ref{theorem_result_dis_expectation}, we obtain the results.
	
\end{proof}

\begin{proof}[\bf Proof of Corollary \ref{cor_consistency_in_ell2}]
	
From Lemma 	\ref{lm_b1_Xie_lwerbound_hell}, we have for some positive constant $ K $ that
	\begin{align*}
K \left[
\mathbb{E}_{\rm x} [ X_1^\top(\theta-\theta_0 )]^2
+ (\sigma- \sigma_0)^2 
\right]
\leq
H^2(P_{\theta,\sigma} ,P_{\theta_0,\sigma_0} )
.
\end{align*}
Thus as a consequence from Corollary \ref{cor_Hellinger_consistency}, we get that
	\begin{align*}
\mathbb{E}
\left[
\int
K (\mathbb{E}_{\rm x} [X_1^\top(\theta-\theta_0 )]^2
+
(\sigma- \sigma_0)^2 
) 
\Pi_{n,\alpha}({\rm d}\theta , {\rm d} \sigma )  
\right]
\leq 
\kappa_\alpha 
\varepsilon_n 
.
\end{align*}
By observing that
$
a^2 \leq a^2 + b^2
$, one deduces that
	\begin{align*}
\mathbb{E}\left[\int
 \mathbb{E}_{\rm x} [X_1^\top(\theta-\theta_0 )]^2
\Pi_{n,\alpha}({\rm d}\theta , {\rm d} \sigma )  
\right]
\leq 
 \frac{\kappa_\alpha}{K}
\varepsilon_n 
,
\quad
\quad
\mathbb{E}\left[\int
 (\sigma- \sigma_0)^2 
\Pi_{n,\alpha}({\rm d}\theta , {\rm d} \sigma )  
\right]
\leq 
 \frac{\kappa_\alpha}{K}
\varepsilon_n 
.
\end{align*}
Under Assumption \ref{asume_minimaleigenvalue_estimebound},
we have that
$ 
\lambda_{min}^G \| \theta-\theta_0 \|_2^2
\leq
\mathbb{E}_{\rm x} [X_1^\top(\theta-\theta_0 )^2]
 $
thus, it yields that
	\begin{align*}
\mathbb{E}
\left[
\int
\| \theta-\theta_0 \|_2^2	
\Pi_{n,\alpha}({\rm d}\theta , {\rm d} \sigma )  
\right]
\leq 
\frac{\kappa_\alpha}{\lambda_{min}^G K}
\varepsilon_n 
.
\end{align*}
The proof is completed.
\end{proof}

\subsubsection{Proofs for Subsection \ref{sc_rslton_distribu}}

\begin{proof}[\bf Proof of Theorem~\ref{theorem_main}]
	\label{proof_theorem_main}	

We will utilize the general result from Theorem 2.4 in \cite{alquier2020concentration} and thus must verify its requirements.
	
First, from Lemma \ref{lm_boundfor_variance_KL}, we have that
\begin{align*}
{\rm Var}_{\theta_0,\sigma_0} \left[\log\frac{p_{\theta_0,\sigma_0} }{p_{\theta,\sigma} } (Z_1) \right]
	= 
\frac{(\sigma^2 -\sigma_0^2)^2 }{ 2\sigma^4}
+
\frac{\sigma_0^2 \, \mathbb{E}_{\rm x} (X_1^\top(\theta-\theta_0 )^2) }{\sigma^4} 
 \leq
\frac{(\sigma^2 -\sigma_0^2)^2 }{ 2\sigma^4}
+
 \frac{\sigma_0^2 \, C_{\rm x}^2 
 	\| \theta-\theta_0\|_2^2  }{\sigma^4} 
,
\end{align*}
where we have used Assumption \ref{asume_finite_2ndmoment}.
When integrating with respect to $\rho_n := p_0\otimes \tilde{\pi}_\sigma $, given in \eqref{eq_specific_distribution} and \eqref{eq_prior_sigma_specific}, and using Lemma \ref{lema_boundfor_ell2} and Lemma \ref{lm_bound_for_kl_term}, we have that
	\begin{align*}
	\int
{\rm Var}_{\theta_0,\sigma_0} \left[\log\frac{p_{\theta_0,\sigma_0} }{p_{\theta,\sigma} } (Z_1) \right]
	\rho_n({\rm d} \theta)
	\leq 
\frac{\eta^2 }{ 2\sigma_{min}^{-2} / 16} 
	+
   \frac{\sigma_0^2 \, C_{\rm x}^2 
  	  4\tau^2 d }{\sigma_{min}^{-2} / 16 } 
  .
	\end{align*}
	Moreover, from the proof of Theorem \ref{theorem_result_dis_expectation}, we have that
	\begin{align*}
	\int\mathcal{K}(P_{\theta_0,\sigma_0} ,P_{\theta,\sigma} )\rho_n({\rm d} \theta)
	\leq 
\frac{2C_{\rm x}^2 }{\sigma_{min}^{-1}} 
4d \tau^2 
+ K  
\frac{\log(\varrho) \eta^2}{ \sigma_{min}^{-2} }	.
	\end{align*}	
From Lemma \ref{lema_boundof_KL} and Lemma \ref{lm_sigma_prior}, we have that 
\begin{align*}
\frac1n \mathcal{K}(\rho_n,\pi)
 =
\frac1n \mathcal{K}(p_0,\pi_{\theta} )
+ 
\frac1n \mathcal{K} (\tilde{\pi}_\sigma, \pi_\sigma)
\leq
\frac{	4 s^* \log \left(\frac{C_1 }{\tau s^*}\right)
	+	\log(2) 
}{n}
+
\frac{
	\log \left( \frac{e 2^{a+1}\Gamma(a) }{ \eta^a} \right)
}{n}
.
\end{align*}

Setting $ \tau = \frac{1}{\sqrt{nd} C_{\rm x}} $ and $ \eta^2 = \frac{1}{n} $ yields the following.
	\begin{align*}
\int
{\rm Var}_{\theta_0,\sigma_0} \left[\log\frac{p_{\theta_0,\sigma_0} }{p_{\theta,\sigma} } (Z_1) \right]
\rho_n({\rm d} \theta)
\leq 
 \frac{8}{n\sigma_{min}^{-2} }
+
\frac{ 32\sigma_{max}^{-1}}{n\sigma_{min}^{-2}}
,
\end{align*}
\begin{align*}
\int\mathcal{K}(P_{\theta_0,\sigma_0} ,P_{\theta,\sigma} )\rho_n({\rm d}\theta , {\rm d} \sigma )
	\leq
\frac{8 }{n\sigma_{min}^{-1}}
+ K  
\frac{\log(\varrho)}{n \sigma_{min}^{-2} }
,
\end{align*}
and
\begin{align*}
\frac1n \mathcal{K}(\rho_n,\pi)
\leq
\frac{	4 s^* \log \left(\frac{C_{\rm x} C_1 \sqrt{nd} }{ s^*}\right)
	+	\log(2)
	+ \log \left( e 2^{a+1}\Gamma(a) n^{a/2} \right)}{n}
.
\end{align*}
	Consequently, we can now apply Theorem 2.4 from \cite{alquier2020concentration},  with
	\[
	\varepsilon_n
	=
	c
	\frac{ s^* \log 
		\left(\frac{ C_{\rm x} C_1 \sqrt{nd} }{ s^*}\right) }{n}
	,
	\]
	for some constant $ c>0 $ depending only on $  \sigma_{min},\sigma_{max}, a $.
	The proof is completed.
\end{proof}

\begin{proof}[\bf Proof of Corollary \ref{cor_concentration_Hellinger}]
	From \cite{van2014renyi}, we have that 
	$$ 
	H^2(P_{\theta,\sigma} ,P_{\theta_0,\sigma_0} )
	\leq 
	D_{1/2}(P_{\theta,\sigma} ,P_{\theta_0,\sigma_0} )  
	\leq 
	D_{\alpha}(P_{\theta,\sigma} ,P_{\theta_0,\sigma_0} ) 
	,
	$$
	for $ \alpha \in [0.5,1) $. In addition, 	for $ \alpha \in (0, 0.5)  $, we also have that 
	$$ 
	D_{1/2}(P_{\theta,\sigma} ,P_{\theta_0,\sigma_0} )  
	\leq 
	\frac{(1-\alpha)1/2}{\alpha (1-1/2)} D_{\alpha}(P_{\theta,\sigma} ,P_{\theta_0,\sigma_0} ) 
	=
	\frac{(1-\alpha)}{\alpha} D_{\alpha}(P_{\theta,\sigma} ,P_{\theta_0,\sigma_0} ) 
	.
	$$ 
	
	Thus, using definition of $ \kappa_\alpha  $ and Corollary \ref{cor_concentration}, we obtain the results.	
\end{proof}

\begin{proof}[\bf Proof of Corollary \ref{cor_concen_rate_in_ell2}]
	
The proof follows a similar approach to the proof of Corollary \ref{cor_consistency_in_ell2}. Specifically, we utilize Lemma \ref{lm_b1_Xie_lwerbound_hell} within Corollary \ref{cor_concentration_Hellinger}, leading to the desired results straightforwardly.
\end{proof}

\subsubsection{Proofs for Subsection \ref{sc_mispecifiedcase}}

\begin{proof}[\bf Proof of Theorem~\ref{theorem_misspecified}]

We proceed to apply Theorem 2.7 from \cite{alquier2020concentration} and will now verify its conditions.

First, we have, for $ u \sim p(\theta_0, \sigma_0) $ and with $ \mu = X^\top \theta , \mu_* = X^\top \theta^* , \mu_0 = X^\top \theta_0 $, that
\begin{align*}
	\mathbb{E}_{\theta_0, \sigma_0}\left[
\log\frac{ p_{\theta^*,\sigma_0}}{
	p_{\theta,\sigma_0} } (Z_1)
\right]
& \leq
	\mathbb{E}_{\theta_0, \sigma_0}
	\left[
\frac{ (u-\mu)^2 - (u-\mu_*)^2 }{2\sigma_0^2}
\right]
\\
& \leq
\mathbb{E}_{\theta_0, \sigma_0}
\left[
\frac{ (\mu_*-\mu)(2u-\mu_*- \mu) }{2\sigma_0^2}
\right]
\\
& \leq
\mathbb{E}_{\rm x}
\left[
\frac{ (\mu_*-\mu)(2\mu_0 -\mu_*- \mu) }{2\sigma_0^2}
\right]
\\
& \leq
\mathbb{E}_{\rm x}
\left[
\frac{ (\mu_*-\mu)^2 + 2(\mu_*-\mu)(\mu_0 -\mu_*) + (\mu_0 -\mu_*)^2 }{2\sigma_0^2}
\right]
\\
& \leq
\frac{ \mathbb{E}_{\rm x} (\mu_0 -\mu)^2 }{2\sigma_0^2}
 \leq
\sigma_{min} 
\mathbb{E}_{\rm x}[ X_1^\top (\theta_0 -\theta)]^2 
 \leq
\sigma_{min}  C_{\rm x}^2\|\theta_0 -\theta\|^2 
,
\end{align*}
where we have used Assumption \ref{asume_finite_2ndmoment}  in the last inequality.
	
When integrating with respect to $\rho_n = p_0\otimes \tilde{\pi}_\sigma $, as defined in equations \eqref{eq_specific_distribution} and \eqref{eq_prior_sigma_specific}, and employing Lemma \ref{lema_boundfor_ell2}, we obtain the following:
\begin{align*}
\int \mathbb{E}_{\theta_0, \sigma_0}\left[
\log\frac{ p_{\theta^*,\sigma_0}}{
	p_{\theta,\sigma_0} } (Z_1)\right] \rho_n({\rm d}\theta , {\rm d} \sigma )
	\leq 
\sigma_{min} 
C_{\rm x}^2 4d \tau^2 
.
\end{align*} 
From Lemma \ref{lema_boundof_KL} and Lemma \ref{lm_sigma_prior}, we have that 
\begin{align*}
\frac1n \mathcal{K}(\rho_n,\pi)
=
\frac1n \mathcal{K}(p_0,\pi_{\theta} )
+ 
\frac1n \mathcal{K} (\tilde{\pi}_\sigma, \pi_\sigma)
\leq
\frac{	4 s^* \log \left(\frac{C_1 }{\tau s^*}\right)
	+	\log(2)
}{n}
+ 
\frac{
	\log \left( \frac{e 2^{a+1}\Gamma(a) }{ \eta^a} \right)
}{n}
.
\end{align*}
Putting $ \tau = \frac1{\sqrt{nd} C_{\rm x}}\, $ and $ \eta^2 = 1/n $, it leads to
\begin{align*}
\int \mathbb{E}_{\theta_0, \sigma_0}
\left[
\log\frac{ p_{\theta^*,\sigma_0}}{
	p_{\theta,\sigma_0} } (Z_1)\right]
\rho_n({\rm d}\theta , {\rm d} \sigma )
&	\leq
\frac{4 \sigma_{min} }{n}
,
\end{align*}
and
\begin{align*}
\frac1n \mathcal{K}(\rho_n,\pi)
\leq
\frac{	4 s^* \log \left(\frac{C_{\rm x} C_1 \sqrt{nd} }{ s^*}\right)
	+	\log(2)
	+ \log \left( e 2^{a+1}\Gamma(a) n^{a/2} \right)}{n}
.
\end{align*}
	To obtain an estimate of the rate $\varepsilon_n$ as in Theorem 2.7 in \cite{alquier2020concentration}, we put together those bounds and apply 
	the theorem with
	\[
	\varepsilon_n
	=
c	s^* \log \left(\frac{C_{\rm x} C_1 \sqrt{nd} }{ s^*}\right) / n
	.
	\]
	for some constant $ c>0 $ depending only on $ \sigma_{min} , a $.
The proof is completed.
\end{proof}

\begin{proof}[\bf Proof of Corollary \ref{cor_misspecified}]
	From Lemma 	\ref{lm_b1_Xie_lwerbound_hell}, we have for some positive constant $ K $ that
	\begin{align*}
K \mathbb{E}_{\rm x} [X_1^\top(\theta-\theta_0 )]^2
	\leq
	K \left[
	\mathbb{E}_{\rm x} [X_1^\top(\theta-\theta_0 )]^2
	+ (\sigma- \sigma_0)^2 
	\right]
	\leq
	H^2(P_{\theta,\sigma} ,P_{\theta_0,\sigma_0} )
	.
	\end{align*}
	and from the proof of Corollary \ref{cor_Hellinger_consistency}, there exist some constant $ K_\alpha >0 $ that $ 	H^2(P_{\theta,\sigma} ,P_{\theta_0,\sigma_0} ) 
	\leq
	K_\alpha D_{\alpha} (P_{\theta,\sigma} ,P_{\theta_0,\sigma_0} ) $. Moreover, one has that
	\begin{align*}
	\mathcal{K}(P_{\theta_0,\sigma_0} ,P_{\theta,\sigma_0} )
	=
	\frac{\mathbb{E}_{\rm x} [X_1^\top(\theta-\theta_0 )]^2 }{2\sigma_0^2}
	\leq 
	\sigma_{min} \mathbb{E}_{\rm x} [X_1^\top(\theta-\theta_0 )]^2
	.
	\end{align*}
	Putting these bounds into Theorem \ref{theorem_misspecified}, it yields that
	\begin{equation*}
	\mathbb{E} \left[ \int
	\mathbb{E}_{\rm x} [X_1^\top(\theta-\theta_0 )]^2
	\Pi_{n,\alpha}({\rm d}\theta , {\rm d} \sigma )\right]
	\leq 
	K_\alpha 
	\left\{
	\min_{\theta\in\Theta} \sigma_{min} \mathbb{E}_{\rm x} [X_1^\top(\theta-\theta_0 )]^2
	+  \varepsilon_n
	\right\}
	,
	\end{equation*}
	for some positive constant $ K_\alpha  $ depending only on $ \alpha $. The proof is completed.
\end{proof}

\subsection{Proof of Section \ref{sc_regular_poster_n_horseshoe}}

\subsubsection{Proof of Subsection \ref{sc_regular_posterior}}

\begin{proof}[\bf Proof of Theorem~\ref{thm_regular_posterior}]
	\label{proof_thm_regular_posterior}	
	
Based on the relationship described in \eqref{eq_connection_fracpost_regu_post} between the regular posterior and the fractional posterior, we will establish the concentration characteristic of the fractional posterior $\Pi_{n,\alpha}( \theta, \sigma_* )$, which will consequently imply the corresponding property of the regular posterior.
	We can check the hypotheses on the KL between the likelihood terms as required in Theorem 2.4 in \cite{alquier2020concentration}. 
	
	From the proof of Theorem \ref{theorem_main}, with $\rho_n := p_0\otimes \tilde{\pi}_{\sigma_*} $, given in \eqref{eq_specific_distribution} and \eqref{eq_prior_sigmastar_specific}, by using Lemma \ref{lema_boundfor_ell2} and Lemma \ref{lm_bound_for_kl_term}, we have that
	\begin{align*}
	\int
	{\rm Var}_{\theta_0,\sigma_0} \left[\log\frac{p_{\theta_0,\sigma_0} }{p_{\theta,\sigma} } (Z_1) \right]
	\rho_n({\rm d} \theta)
	\leq 
\frac{\eta^2 }{ 2\sigma_{min}^{-2} / 16} 
+
\frac{\sigma_0^2 \, C_{\rm x}^2 
	4\tau^2 d }{\sigma_{min}^{-2} / 16 } 
	.
	\end{align*}
	Moreover, from the proof of Theorem \ref{theorem_result_dis_expectation}, we have that
	\begin{align*}
	\int\mathcal{K}(P_{\theta_0,\sigma_0} ,P_{\theta,\sigma} )\rho_n({\rm d} \theta)
	\leq 
\frac{2C_{\rm x}^2 }{\sigma_{min}^{-1}} 
4d \tau^2 
+ K  
\frac{\log(\varrho) \eta^2}{ \sigma_{min}^{-2} }
	.
	\end{align*}	
Put $ \alpha = 1- 1/(\log n)^t $ for  $ t>1 $. From Lemma \ref{lema_boundof_KL} and Lemma \ref{lm_sigma_prior_regular_poster}, we have that 
	\begin{align*}
	\frac1n \mathcal{K}(\rho_n,\pi)
 =
	\frac1n \mathcal{K}(p_0,\pi_{\theta} )
	+ 
	\frac1n \mathcal{K} (\tilde{\pi}_{\sigma_*}, \pi_{\sigma_*})
	\leq
	\frac{	4 s^* \log \left(\frac{C_1 }{\tau s^*}\right)
		+	\log(2) + 
K n \eta
	}{n}
	.
	\end{align*}
	
	Putting $ \tau = \frac1{\sqrt{nd} C_{\rm x}}\, $ and $ \eta = 1/n $, it leads to
	\begin{align*}
	\int
	{\rm Var}_{\theta_0,\sigma_0} \left[\log\frac{p_{\theta_0,\sigma_0} }{p_{\theta,\sigma} } (Z_1) \right]
	\rho_n({\rm d} \theta)
	\leq 
 \frac{8}{n^2 \sigma_{min}^{-2} }
+
\frac{ 32\sigma_{max}^{-1}}{n\sigma_{min}^{-2}}
	,
	\end{align*}
	\begin{align*}
	\int\mathcal{K}(P_{\theta_0,\sigma_0} ,P_{\theta,\sigma} )\rho_n({\rm d}\theta , {\rm d} \sigma )
		\leq
\frac{8 }{n\sigma_{min}^{-1}}
+ K  
\frac{\log(\varrho)}{n^2 \sigma_{min}^{-2} }
	,
	\end{align*}
	and
	\begin{align*}
	\frac1n \mathcal{K}(\rho_n,\pi)
	\leq
	\frac{	4 s^* \log \left(\frac{C_{\rm x} C_1 \sqrt{nd} }{ s^*}\right)
		+	\log(2)
		+ K }{n}
	.
	\end{align*}
	Consequently, we can now apply Theorem 2.6 and Theorem 2.4 from \cite{alquier2020concentration}, with $\rho_n := p_0 $ and  with
	\[
	\varepsilon_n
	=
	c s^* \log 	\left(\frac{ C_{\rm x} C_1 \sqrt{nd} }{ s^*}\right) / n
	,
	\]
	for some constant $ c>0 $ depending only on $  \sigma_{min},\sigma_{max}, a $. Thus, we obtain concentration result for the fractional posterior $ \Pi_{n,\alpha}( \theta, \sigma_* ) $. This means that the concentration result for the regular posterior $ \Pi_{n}(\theta,\sigma ) $ is obtained.
	The proof is completed.	
\end{proof}

\begin{proof}[\bf Proof of Theorem \ref{thm_misspecified_regularposte}]
	
	The proof follows a similar argument as in the proof of Theorem \ref{theorem_misspecified}. We proceed to apply Theorem 2.7 from \cite{alquier2020concentration} and will now verify its conditions. From the proof of Theorem \ref{theorem_misspecified}, we have that
	\begin{align*}
	\mathbb{E}_{\theta_0, \sigma_0}\left[
	\log\frac{ p_{\theta^*,\sigma_0}}{
		p_{\theta,\sigma_0} } (Z_1)
	\right]
	\leq
	\sigma_{min}  C_{\rm x}^2\|\theta_0 -\theta\|^2 
	.
	\end{align*}

	When integrating with respect to $\rho_n = p_0\otimes \tilde{\pi}_{\sigma_*} $, as defined in equations \eqref{eq_specific_distribution} and \eqref{eq_prior_sigmastar_specific}, and employing Lemma \ref{lema_boundfor_ell2}, we obtain the following:
	\begin{align*}
	\int \mathbb{E}_{\theta_0, \sigma_0}\left[
	\log\frac{ p_{\theta^*,\sigma_0}}{
		p_{\theta,\sigma_0} } (Z_1)\right] \rho_n({\rm d}\theta , {\rm d} \sigma )
	\leq 
	\sigma_{min} 
	C_{\rm x}^2 4d \tau^2 
	.
	\end{align*} 
Put $ \alpha = 1- 1/(\log n)^t $ for  $ t>1 $. From Lemma \ref{lema_boundof_KL} and Lemma \ref{lm_sigma_prior_regular_poster}, we have that 
\begin{align*}
\frac1n \mathcal{K}(\rho_n,\pi)
=
\frac1n \mathcal{K}(p_0,\pi_{\theta} )
+ 
\frac1n \mathcal{K} (\tilde{\pi}_{\sigma_*}, \pi_{\sigma_*})
\leq
\frac{	4 s^* \log \left(\frac{C_1 }{\tau s^*}\right)
	+	\log(2) + 
	K n \eta
}{n}
.
\end{align*}
	Putting $ \tau = \frac1{\sqrt{nd} C_{\rm x}}\, $ and $ \eta = 1/n $, it leads to
	\begin{align*}
	\int \mathbb{E}_{\theta_0, \sigma_0}
	\left[
	\log\frac{ p_{\theta^*,\sigma_0}}{
		p_{\theta,\sigma_0} } (Z_1)\right]
	\rho_n({\rm d}\theta , {\rm d} \sigma )
	&	\leq
	\frac{4 \sigma_{min} }{n}
	,
	\end{align*}
	\begin{align*}
\frac1n \mathcal{K}(\rho_n,\pi)
\leq
\frac{	4 s^* \log \left(\frac{C_{\rm x} C_1 \sqrt{nd} }{ s^*}\right)
	+	\log(2)
	+ K }{n}
.
\end{align*}
	To obtain an estimate of the rate $\varepsilon_n$ as in Theorem 2.7 in \cite{alquier2020concentration}, we put together those bounds and apply 
	the theorem with
	\[
	\varepsilon_n
	=
	c	s^* \log \left(\frac{C_{\rm x} C_1 \sqrt{nd} }{ s^*}\right) / n
	.
	\]
	for some constant $ c>0 $ depending only on $ \sigma_{min} , a $.
	The proof is completed.	
\end{proof}

\subsubsection{Proofs for Subsection \ref{sc_spike_n_slab_prior}}

\begin{proof}[\bf Proof of Theorem \ref{theorem_spike_slab}] 
	
Under the assumption that $\|\theta_0\|_2 = 1 $.
We start by defining, for any  $ \eta > 0 $, 
	$$ 
	\rho_{\theta_0,\eta}
	({\rm d}\theta) 
	\propto
\mathbf{1}_{\|\theta-\theta_0\|_2 \leq \eta}\, 
	\pi_{SnS} ({\rm d}\theta) 
	.
	$$
First, from Lemma \ref{lm_boundfor_variance_KL}, we have that
\begin{align*}
{\rm Var}_{\theta_0,\sigma_0} \left[\log\frac{p_{\theta_0,\sigma_0} }{p_{\theta,\sigma} } (Z_1) \right]
 \leq
\frac{(\sigma^2 -\sigma_0^2)^2 }{ 2\sigma^4}
+
\frac{\sigma_0^2 \, C_{\rm x}^2 
	\| \theta-\theta_0\|_2^2  }{\sigma^4} 
,
\end{align*}
where we have used Assumption \ref{asume_finite_2ndmoment}.
When integrating with respect to $\rho_n := 	\rho_{\theta_0,\eta} \otimes \tilde{\pi}_\sigma $, where $ \tilde{\pi}_\sigma $ given in \eqref{eq_prior_sigma_specific}, and using Lemma \ref{lm_bound_for_kl_term}, we have that
\begin{align*}
\int
{\rm Var}_{\theta_0,\sigma_0} \left[\log\frac{p_{\theta_0,\sigma_0} }{p_{\theta,\sigma} } (Z_1) \right]
\rho_n({\rm d} \theta)
\leq 
\frac{\eta^2 }{ 2\sigma_{min}^{-2} / 16} 
+
\frac{\sigma_0^2 \, C_{\rm x}^2 
\eta^2 }{\sigma_{min}^{-2} / 16 } 
.
\end{align*}
Moreover, from the proof of Theorem \ref{theorem_result_dis_expectation}, we have that
\begin{align*}
\int\mathcal{K}(P_{\theta_0,\sigma_0} ,P_{\theta,\sigma} )\rho_n({\rm d} \theta)
\leq 
\frac{2C_{\rm x}^2 }{\sigma_{min}^{-1}} 
\eta^2
+ K  
\frac{\log(\varrho) \eta^2}{ \sigma_{min}^{-2} }	.
\end{align*}	
From Lemma \ref{lema_boundof_KL} and Lemma \ref{lm_sigma_prior}, we have that 
\begin{align*}
\frac1n \mathcal{K}(\rho_n,\pi)
 =
\frac1n \mathcal{K}(p_0,\pi_{\theta} )
+ 
\frac1n \mathcal{K} (\tilde{\pi}_\sigma, \pi_\sigma)
\leq
\frac{	K_{v_1} s^* \log 
	\left(
	\frac{ \sqrt{d} }{\eta}
	\right)	+ 1 }{n}
+
\frac{	\log \left( \frac{e 2^{a+1}\Gamma(a) 
	}{ \eta^a} \right)
}{n}
.
\end{align*}

Setting $ \eta^2 = \frac{s^*}{n C_{\rm x}^2 } $ yields the following.
\begin{align*}
\int
{\rm Var}_{\theta_0,\sigma_0} \left[\log\frac{p_{\theta_0,\sigma_0} }{p_{\theta,\sigma} } (Z_1) \right]
\rho_n({\rm d} \theta)
\leq 
\frac{8s^*}{n\sigma_{min}^{-2} C_{\rm x}^2 }
+
\frac{ 32\sigma_{max}^{-1}s^* }{n\sigma_{min}^{-2}}
,
\end{align*}
\begin{align*}
\int\mathcal{K}(P_{\theta_0,\sigma_0} ,P_{\theta,\sigma} )\rho_n({\rm d}\theta , {\rm d} \sigma )
\leq
\frac{2 s^* }{n\sigma_{min}^{-1}}
+ K  
\frac{\log(\varrho) s^* }{n \sigma_{min}^{-2} C_{\rm x}^2 }
,
\end{align*}
and
\begin{align*}
\frac1n \mathcal{K}(\rho_n,\pi)
\leq
\frac{	4 s^* \log \left(\frac{C_{\rm x} \sqrt{nd} }{ s^*}\right)	+ 1
	+ \log \left( e 2^{a+1}\Gamma(a) n^{a/2} C_{\rm x}^a/s^* \right)}{n}
.
\end{align*}
Consequently, we can now apply Theorem 2.4 from \cite{alquier2020concentration},  with
\[
\varepsilon_n
=
c s^* \log 	\left(\frac{ C_{\rm x} \sqrt{nd} }{ s^*}\right) / n
,
\]
for some constant $ c>0 $ depending only on $  \sigma_{min},\sigma_{max}, a $.
The proof is completed.
\end{proof}

\subsection{Useful Lemmas}
%%% the translate prior
\begin{dfn}
	Let's define the following distribution as a translation of the prior $ \pi_\theta $,
	\begin{equation}
	\label{eq_specific_distribution}
	p_0(\theta) 
	\propto 
	\pi_\theta (\theta - \theta_0 )
\mathbf{1}_{B_1(2d\tau)} (\theta - \theta_0 ).
	\end{equation}
\end{dfn}
Given that \( \| \theta_0 \|_1 \leq C_1 - 2d\tau \) and considering \( \theta - \theta_0 \in B_1(2d\tau) \), it follows that \( \theta \in B_1(C_1) \). Consequently, the distribution \( p_0 \) is absolutely continuous with respect to the prior distribution \( \pi_\theta \), resulting in a finite Kullback-Leibler divergence.

The following two lemmas can be found in \cite{mai2023high}. Their proofs are derived by utilizing Lemma 2 and Lemma 3 from \cite{dalalyan2012mirror}.

\begin{lemma}
	\label{lema_boundof_KL}
	For $p_0$ defined in (\ref{eq_specific_distribution}), we have that
	$$
\mathcal{K}(p_0,\pi_{\theta} )
	\leq
	4 s^* \log \left(\frac{C_1 }{\tau s^*}\right)
	+
	\log(2)
	.
	$$
\end{lemma}
\begin{lemma}
	\label{lema_boundfor_ell2}
	For $p_0$ defined in (\ref{eq_specific_distribution}) and for
	$d\geq 2$, we have that
	$$
	\int \| \theta - \theta_0 \|^2 p_0({\rm d }\theta)
	\leq
	4d\tau^2 
	.
	$$
\end{lemma}

\begin{lemma}
	\label{lm_bound_error_variance_Pati2014}
	Consider the function $ h_{\beta}(x) = \log x - (x-1) + \beta (x-1)^2/2 $ on $(0, \infty)$ for $\beta > 1$. Let $\varepsilon \in (0, 1/2)$ and set $\beta = 8 \log(1/\varepsilon)$. Then $h_{\beta}(x) \geq 0$ for all $x \geq \varepsilon$. 
\end{lemma}
\begin{proof}
	The proof can be found in the proof of Lemma 1.3, in the supplementary document, of \cite{pati2014posterior}. 
\end{proof}

\begin{lemma}
	\label{lm_bound_for_kl_term} 

Under Assumption \ref{asume_sigma_lowerbound} and assuming that \( | \sigma_0^2 - \sigma^2 | \leq \eta \) and \( \eta / ( 2 \sigma_{min})^{-1} < 1/2 \), we have that
	\begin{align*}
 \sigma^2 \geq  \sigma_{min}^{-1} / 4 ,
 \quad 
 \text{ and }	
	\quad
	\log\left( \frac{\sigma^2}{\sigma_0^2} \right) 
	+
	\frac{\sigma_0^2}{\sigma^2} - 1
	\leq
	K  
	\frac{\log(\varrho) }{\sigma_{min}^{-2} 
	}\eta^2 
	\end{align*}
	for some absolute constant $ K>0 $ and $ \varrho = 2  \sigma_{max}^{-1}/\sigma_{min}^{-1}  $.
\end{lemma}

\begin{proof}[\bf Proof of Lemma \ref{lm_bound_for_kl_term}]
	Applying Lemma \ref{lm_bound_error_variance_Pati2014}, for $ x = \frac{\sigma_0^2}{\sigma^2} $, with $\varepsilon = \frac{( 2 \sigma_{min})^{-1} }{2 ( 2 \sigma_{max})^{-1}} < 1/2  $ and $\beta = 8 \log(1/\varepsilon) $, one gets that
	\begin{align*}
	\log\left( \frac{\sigma_0^2}{\sigma^2} \right) 
	-
	\left( 
	\frac{\sigma_0^2}{\sigma^2} - 1
	\right)
	\geq
	- 4 \log(1/\varepsilon) 	\left(\frac{\sigma_0^2}{\sigma^2} -1\right)^2
	.
	\end{align*}
Under Assumption \ref{asume_sigma_lowerbound},
	$
	| \sigma_0^2 - \sigma^2 | 
	\geq
	\sigma_0^2
	\left|
	\frac{\sigma^2}{\sigma_0^2} -1
	\right|
	\geq
	( 2 \sigma_{min})^{-1}
	\left|
	\frac{\sigma^2}{\sigma_0^2} -1
	\right|
	,
	$
	and thus 
$
	\left|
	\frac{\sigma^2}{\sigma_0^2} -1
	\right|
	\leq
	\eta / ( 2 \sigma_{min})^{-1}
	.
$
As a consequence,
	$
	\frac{\sigma^2}{\sigma_0^2}
	\geq
	1
	-\eta / ( 2 \sigma_{min})^{-1}
	$, and then from Assumption \ref{asume_sigma_lowerbound},
	\begin{align*}
	\sigma^2
	\geq
	( 2 \sigma_{min})^{-1}
	(	1 -\eta / ( 2 \sigma_{min})^{-1}) 
	\geq
	( 2 \sigma_{min})^{-1}/2
	.
	\end{align*}
	Putting all together, 
	\begin{align*}
	\left(
	\frac{\sigma_0^2}{\sigma^2}-1 \right)^2
	=
	\left(
	\frac{\sigma_0^2 - \sigma^2}{\sigma^2} 
	\right)^2
	\leq
	\left(
	\frac{\eta}{ ( 2 \sigma_{min})^{-1}/2 } 
	\right)^2
	.	
	\end{align*}
	Thus, using Lemma \ref{lm_bound_error_variance_Pati2014}, there exists a constant $ K>0 $ such that
	\begin{align*}
	\log\left( \frac{\sigma_0^2}{\sigma^2} \right) 
	-
	\left( 
	\frac{\sigma_0^2}{\sigma^2} - 1
	\right)
	\geq
	- K \log(1/\varepsilon) 
	\frac{\eta^2}{	\left( ( 2 \sigma_{min})^{-1}  
		\right)^2}
	\end{align*}
	we obtain the result. The proof is completed.	
\end{proof}

\begin{lemma}
	\label{lm_sigma_prior}
	Let $ \sigma^2 \sim IG(a, b)$ for some fixed $ a>0 $. Put 
	\begin{align}
	\label{eq_prior_sigma_specific}
	\tilde{\pi}_\sigma 
	\propto 
	\mathbf{1}_{\vert \sigma^2 - \sigma_0^2 \vert < \eta}\, \pi_\sigma
	.
	\end{align}
	 Then for any fixed $\sigma_0^2 > 0$ and $ \eta >0 $, we have for $ b=\eta $ that
	\begin{align*}
	\mathcal{K} (\tilde{\pi}_\sigma , \pi_\sigma )
	\leq
	\log( \frac{e 2^{a+1}\Gamma(a) }{ \eta^a} )	
	\end{align*}
\end{lemma}
\begin{proof}[\bf Proof of Lemma \ref{lm_sigma_prior}]
	One has that
	\begin{align*}
	\mathcal{K} (\tilde{\pi}_\sigma , \pi_\sigma )
	=
	\log\frac{1}{ \pi_\sigma ( \vert \sigma^2 - \sigma_0^2 \vert < \eta) },
	\end{align*}	
	thus to upper bound the KL, we just need to lower bound $ \pi_\sigma ( \vert \sigma^2 - \sigma_0^2 \vert < \eta) $.
	
	Without loss of generality let $\sigma_0^2 = 1$. Since otherwise 
	$ \pi_\sigma [\,\,\vert \sigma^2 - \sigma_0^2 \vert < \eta] 
	= 
	\pi_\sigma [\,\,\vert \sigma^2/\sigma_0^2 - 1 \vert < \eta/\sigma_0^2] 
	= 
	\pi_\sigma [\,\,\vert \tilde\sigma^2 - 1 \vert < \delta],$ where $ \tilde\sigma^2 \sim IG( a, b/\sigma_0^2 )$ and $\delta = \eta/\sigma_0^2$ is fixed.
	\begin{align*}
	\pi_\sigma (\,\vert \sigma^2 - 1 \vert < \eta) 
	&
	\geq 
	\pi_\sigma (1 - \eta < \sigma^2 < 1+\eta)
\geq 
	\pi_\sigma ( \eta < \sigma^2 < 2\eta)
	.
	\end{align*}
	and
	\begin{align*}
\pi_\sigma ( \eta < \sigma^2 < 2\eta)
& =  
\dfrac{b^a}{\Gamma \{a\}}
\int_{\eta}^{2\eta} x^{-a-1}\exp\left(-\frac{b}{ x}\right) {\rm d}x
\\
& \geq 
\dfrac{b^a}{\Gamma \{a\}}
\eta (2 \eta)^{-a-1}
\exp\left(-\frac{b}{\eta}\right)
 =
\frac{b^a}{ 2^{a+1}\Gamma(a)} 
\exp\left[- \frac{b}{\eta}\right] \eta^{-a}
 \geq 
\frac{b^a}{ 2^{a+1}\Gamma(a)} 
\exp\left[- \frac{b}{\eta}\right]
.
\end{align*}	
	Therefore,
	\begin{align*}
	\mathcal{K} (\tilde{\pi}_\sigma , \pi_\sigma )
	\leq
	\log( \frac{ 2^{a+1}\Gamma(a) }{ b^a} )
	+ \frac{b}{\eta}
	,
	\end{align*}
	choosing $ b = \eta $, we obtain the result. The proof is completed.
\end{proof}

\begin{lemma}[Lemma S7 in Supplementary material of \cite{chakraborty2020bayesian}]
	\label{lm_sigma_prior_regular_poster}
	Let $\tau^2 \sim IG(n(1-\alpha)/2 + a, \alpha b)$ for some fixed $a , b >0$ and $\alpha = \{1- 1/(\log n)^t\}, t>1$.  Then, there exist some constant $K>0 $, for any fixed $\sigma_0^2 > 0 $ and $\eta >0$, one has that
	$
	P[\,\,\vert \tau^2 - \sigma_0^2 \vert < \eta] 
	\geq e^{-K n \eta }.$	And,  put 
	\begin{align}
	\label{eq_prior_sigmastar_specific}
	\tilde{\pi}_{\sigma_*} 
	\propto 
	\mathbf{1}_{\vert \sigma^2_* - \sigma_0^2 \vert < \eta}\, \pi_\sigma
	.
	\end{align}
	Then for any fixed $\sigma_0^2 > 0$ and $ \eta >0 $, we have that
$
	\mathcal{K} (\tilde{\pi}_{\sigma_*} , \pi_\sigma )
	\leq
K n \eta 	
$.
\end{lemma}

\begin{lemma}
	\label{lm_boundfor_variance_KL}
	Let $ p $ and $ p_0 $ be the densities of two Gaussian distributions $ \mathcal{N}(\mu,\sigma^2) $ and $ \mathcal{N}(\mu_0,\sigma_0^2) $, respectively. Then,
	\begin{align*}
	{\rm Var}_{p_0} \left( \log \frac{p_0}{p} \right) 
	=
	\frac{\sigma_0^2}{\sigma^4} (\mu_0 -\mu)^2 
	+
	\frac{(\sigma^2 -\sigma_0^2)^2 }{ 2\sigma^4}
	.
	\end{align*}
\end{lemma}

\begin{proof}[\bf Proof of Lemma \ref{lm_boundfor_variance_KL}]
	For $ u\sim \mathcal{N}(\mu_0,\sigma_0^2) $,
	\begin{align*}
	{\rm Var}_{p_0} \left( \log \frac{p_0}{p} \right) 
	& =
	{\rm Var}_{p_0} \left[  
	\frac{(u - \mu)^2 }{2\sigma^2}
	-
	\frac{(u - \mu_0)^2 }{2\sigma_0^2}  
	\right]
	\\
	& =
	{\rm Var}_{p_0} \left[  
	\frac{(u - \mu_0 + \mu_0 - \mu)^2 }{2\sigma^2}
	-
	\frac{(u - \mu_0)^2 }{2\sigma_0^2}  
	\right]
	\\
	& =
	{\rm Var}_{p_0} \left[  
	\frac{1}{2} \left( \frac{\sigma_0^2}{\sigma^2} -1 \right)
	\frac{(u - \mu_0)^2 }{\sigma_0^2} 
	+
	\frac{\sigma_0 (\mu_0 - \mu)}{\sigma^2}
	\frac{(u - \mu_0) }{\sigma_0} 
	\right]
	.
	\end{align*}
	Now, with a change of variable that $ z:=\frac{(u - \mu_0) }{\sigma_0}   $ then $ z \sim  \mathcal{N}(0,1) $, then
	\begin{align*}
	{\rm Var}_{p_0} \left( \log \frac{p_0}{p} \right) 
	 =
	{\rm Var}_{z} \left[  
	\frac{1}{2} \left(\frac{\sigma_0^2}{\sigma^2}-1 \right)
	z^2
	+
	\frac{\sigma_0 (\mu_0 - \mu)}{\sigma^2}
	z
	\right]
 =
	\frac{(\sigma^2 -\sigma_0^2)^2 }{ 2\sigma^4}
	+
	\frac{\sigma_0^2}{\sigma^4} (\mu_0 -\mu)^2 
	.
	\end{align*}
	The proof is completed.
\end{proof}

\begin{lemma}[Lemma B1 in \cite{xie2020adaptive}]
	\label{lm_b1_Xie_lwerbound_hell}
	There exist some universal constant $ K>0 $ such that
	\begin{align*}
	H^2(P_{\theta,\sigma} ,P_{\theta_0,\sigma_0} )
	\geq
	K (\mathbb{E}_{\rm x} [X_1^\top(\theta-\theta_0 )]^2
	+
	 (\sigma- \sigma_0)^2 
	)
	.
	\end{align*}
\end{lemma}

\begin{lemma}
	\label{lm_concentration_spike_slad}
	For any $\theta_0 \in\mathbb{R}^{p}$ such that $\|\theta_0\| \leq 1$ and $\delta \in (0,1) $, 	as soon as $v_0 \leq \delta^2/(2d \log(2))$,  with $ p = 1- 2e^{-1/d} $, there exists a positive constant $K_{v_1}$, dependent solely on $v_1$, such that the following holds.
	\begin{align*}
	\pi_{SnS} \left(\left\{\theta: \|\theta-\theta_0\|^2 \leq \delta^2 \right\}\right) 
	\geq 
	e^{- K_{v_1}
		s^*\log(\sqrt{d}/(p \delta)) }
	.
	\end{align*}	
\end{lemma}

\begin{proof}[\bf Proof of Lemma \ref{lm_concentration_spike_slad}]
	
	\begin{align*}
	\log \pi\left(\left\{\theta: \|\theta-\theta_0\|^2 \leq \delta^2 \right\}\right) 
	& \geq 
	\log \pi\left(\left\{\theta: \forall i, (\theta_i-\theta_{0,i})^2 \leq \frac{\delta^2}{d} \right\}\right) 
	\\
	& = 
	\sum_{i:\theta_{0,i}\neq 0} 
	\log \pi\left((\theta_i-\theta_{0,i})^2 \leq \frac{\delta^2}{d} \right) 
	+ 
	\sum_{i:\theta_{0,i} = 0} 
	\log \pi \left(\left\{
	\theta_i^2 < \frac{\delta^2}{d} \right\}\right) 
	.
	\end{align*}
	Assume first that $i$ is such that
	$\theta_{0,i}= 0$ and put $\pi_0$ denotes the probability distribution such that the $\theta_i$ are iid $\mathcal{N}(0,v_0)$. One has that
	\begin{align*}
	\sum_{i:\theta_{0,i} = 0} 
	\log \pi \left(\left\{
	\theta_i^2 < \frac{\delta^2}{d} \right\}\right) 
	& \geq 
	\sum_{i:\theta_{0,i} = 0} 
	\log \pi_{0} \left(\left\{ 
	\theta_i^2 < \frac{\delta^2}{d}\right\}\right) 
	+ (d-s^*) \log(1-p) 
	\\
	& = 
	\sum_{i:\theta_{0,i} = 0} 
	\log \left[1-\pi_0 \left(\left\{
	\theta_i^2 > \frac{\delta^2}{d} \right\}\right) \right] 
	+ (d-s^*)\log(1-p) 
	.
	\end{align*}	
	As, one has that
	\begin{align*}
	\pi_0 \left( \theta_i^2 > \frac{\delta^2}{d} \right)
	= 
	\pi_0 \left( \left|\frac{\theta_i}{\sqrt{v_0}}\right| > \frac{\delta}{\sqrt{v_0 d}} \right) 
	\leq
	\exp\left( - \frac{\delta^2}{2 v_0 d} \right)
	\leq
	\frac{1}{2}	
	,
	\end{align*}
	as soon as $v_0 \leq \delta^2/(2d \log(2))$. Thus, we obtain that
	\begin{align*}
	\sum_{i:\theta_{0,i} = 0} 
	\log \pi \left(\left\{
	\theta_i^2 < \frac{\delta^2}{d} \right\}\right) 
	& \geq 
	(d-s^*)\log((1-p)/2) 
	.
	\end{align*}

	Now, considering that $i$ is such that
	$\theta_{0,i}\neq 0$ and that $\theta_{0,i}>0$ (the proof is exactly symmetric if $\theta_{0,i}<0$):
	\begin{align*}
	\pi\left(\left\{\theta: (\theta_i-\theta_{0,i})^2 
	\leq 
	\frac{\delta^2}{d}
	\right\}\right) 
	& \geq 
	p \int_{ | \theta_i-\theta_{0,i}| 
		\leq 
		\frac{\delta}{\sqrt{d} }}
	\frac{\exp(-\theta_i^2/(2v_1^2))}{v_1 \sqrt{2\pi}}
	{\rm d} \theta_i 
	,
	\end{align*}
	and
	\begin{align*}
	\int_{ | \theta_i-\theta_{0,i}| 
		\leq 
		\frac{\delta}{\sqrt{d} }}
	\frac{ e^{-\theta_i^2/(2v_1^2)} }{v_1 \sqrt{2\pi}}
	{\rm d} \theta_i 
	\geq 
	\frac{e^{-(1+1)^2/(2v_1^2)} }{v_1 \sqrt{2\pi}}
	\int_{ | \theta_i-\theta_{0,i}| 
		\leq \frac{\delta}{\sqrt{d} }}
	{\rm d} \theta_i 
	\geq 
	\frac{e^{-2/v_1^2} }{v_1 \sqrt{2\pi}}
	\frac{\delta}{\sqrt{d} }
	,
	\end{align*}

	Putting everything together:
	\begin{align*}
	\log \pi\left(\left\{\theta: \|\theta-\theta_0\|^2 \leq \delta^2 \right\}\right) 
	& 
	\geq 
	s^* \log\left(\frac{ p \delta}{2 \sqrt{2 \pi v_1 d}}
	\right) 
	+ (d-s^*)\log((1-p)/2) 
	\\
	& = 
	s^*  \log\left(\frac{ p \delta }{ (1-p) \sqrt{2 \pi v_1 d}}\right) 
	+ 
	d\log((1-p)/2) 
	.
	\end{align*}
	with $ p = 1- 2e^{-1/d} $,
	\begin{align*}
	\log \pi\left(\left\{\theta: \|\theta-\theta_0\|^2 \leq \delta^2 \right\}\right) 
	\geq 
	s^*  \log\left(\frac{ p \delta }{ (1-p) \sqrt{2 \pi v_1 d}}\right) 
	-
	\frac{d}{d}
	\geq 
	K_{v_1} s^*\log\left(
	\frac{ p \delta}{(1-p)\sqrt{d}}
	\right) -1
	,
	\end{align*}	
	for some positive constant $ K_{v_1} $ depending only on $ v_1 $. This completes the proof.	
\end{proof}

\end{document}